\documentclass[a4paper,12pt]{amsart}
\usepackage{url}
\usepackage[margin=0.9in]{geometry}
\usepackage[english]{babel}
\usepackage[utf8]{inputenc}
\usepackage{amsmath}
\usepackage{amssymb}
\usepackage{cite}
\usepackage{graphicx}
\usepackage{booktabs}  
\usepackage{listings}
\usepackage{pgfplots}
\newcommand{\uproman}[1]{\uppercase\expandafter{\romannumeral#1}}
\usepackage{tkz-euclide}
\usetikzlibrary{calc,angles,positioning,intersections,quotes,decorations.markings,shapes}
\usepackage{amsbsy}
\newtheorem{theorem}{Theorem}[section]

\theoremstyle{definition}

\theoremstyle{remark}

\numberwithin{equation}{section}

\begin{document}

\title[Computing Laplacian eigenvalues of 2D shapes with dihedral symmetry]{Computation of Laplacian eigenvalues of two-dimensional shapes with dihedral symmetry}


\author[D. Berghaus]{David Berghaus}
\address{Bethe Center, University of Bonn, Nussallee 12, 53115 Bonn, Germany}
\curraddr{}
\email{berghaus@th.physik.uni-bonn.de}
\thanks{}

\author[R. S. Jones]{Robert Stephen Jones}
\address{Independent Researcher, 43074 Sunbury, Ohio, USA}
\curraddr{}
\email{rsjones7@yahoo.com}
\thanks{}

\author[H. Monien]{Hartmut Monien}
\address{Bethe Center, University of Bonn, Nussallee 12, 53115 Bonn, Germany}
\curraddr{}
\email{hmonien@uni-bonn.de}
\thanks{}

\author[D. Radchenko]{Danylo Radchenko}
\address{Laboratoire Paul Painlevé, University of Lille, F-59655 Villeneuve d'Ascq, France}
\curraddr{}
\email{danradchenko@gmail.com}
\thanks{}


\date{}

\dedicatory{}

\begin{abstract}
We numerically compute the lowest Laplacian eigenvalues of several two-dimensional shapes with dihedral symmetry at arbitrary precision arithmetic. Our approach is based on the method of particular solutions with domain decomposition. We are particularly interested in asymptotic expansions of the eigenvalues $\lambda(n)$ of shapes with $n$ edges that are of the form $\lambda(n) \sim x\sum_{k=0}^{\infty} \frac{C_k(x)}{n^k}$ where $x$ is the limiting eigenvalue for $n\rightarrow \infty$. Expansions of this form have previously only been known for regular polygons with Dirichlet boundary condition and (quite surprisingly) involve Riemann zeta values and single-valued multiple zeta values, which makes them interesting to study. We provide numerical evidence for closed-form expressions of higher order $C_k(x)$ and give more examples of shapes for which such expansions are possible (including regular polygons with Neumann boundary condition, regular star polygons and star shapes with sinusoidal boundary).
\end{abstract}

\maketitle
\allowdisplaybreaks

\section{Introduction}
\label{sec:introduction}
Let $\Psi$ be a function defined on a domain $\Omega\subset\mathbb{R}^2$ that satisfies the Laplace eigenvalue equation
\begin{equation}
\label{eq:laplace_eig}
-\Delta \Psi(x,y) = \lambda \cdot\Psi(x,y)\,,
\end{equation}
where $\Delta = \partial_x^2+\partial_y^2$ denotes the Laplacian in two-dimensional Euclidean space and $\lambda\in\mathbb{R}$. We consider the case when the two-dimensional shape $\Omega$ has the symmetry group of a regular $n$-gon, and $\Psi$ has either Dirichlet or Neumann boundary condition on $\partial \Omega$
\begin{align*}
    \Psi|_{\partial \Omega} &= 0\quad  (\mathrm{Dirichlet})\,,\\
    \partial_{\Vec{n}}\Psi|_{\partial \Omega} &= 0\quad (\mathrm{Neumann})\,,
\end{align*}
where $\partial_{\Vec{n}}$ denotes the normal derivative. By combining Eq.~\eqref{eq:laplace_eig} with the boundary condition, one obtains a discrete spectrum of eigenvalues 
\begin{equation*}
    0<\lambda_1 \leq \lambda_2 \leq \lambda_3 \leq \dots\,,
\end{equation*}
(in the case of Neumann boundary condition we additionally have $\lambda_0=0$). The eigenvalues are invariant under translations and rotations but do depend on the area of the considered shape. More precisely, if $\Omega'$ is obtained from $\Omega$ by a homothety, then one has
\begin{equation}
    \lambda_i(\Omega)A(\Omega) = \lambda_i(\Omega')A(\Omega') \,,
\end{equation}
where $A$ denotes the area of the corresponding shape. Keeping the area constant is therefore crucial for investigating $1/n$ expansions. Throughout this paper, we will always normalize the domains to have area $\pi$, which means that the shapes will approach the unit disk in the limit as $n\to\infty$.

We start by considering regular polygons with Dirichlet boundary condition. Additional examples will be introduced in Section~\ref{sec:further_examples}. There are so far only two regular polygons whose eigenvalues are known explicitly. These are the regular triangle where $\lambda_1 = 4\pi / \sqrt[]{3}$ and the square with $\lambda_1 = 2\pi$, see, e.g., P\'olya and Szeg\"o~\cite{polya} (one should however mention that some eigenmodes of the regular hexagon can be obtained by piecing together eigenmodes of regular triangles). The remaining regular polygons offer challenges, both analytically and numerically, due to the presence of non-analytic vertices (i.e. vertices that are not of the form $\pi/n$). It has been known from works of several authors \cite{Molinari,GRINFELD2012135,boady,jones2017fundamental,analytics_paper} that the lowest eigenvalue of a regular $n$-gon with Dirichlet boundary condition can be approximated by a series in $1/n$
\begin{equation}
\label{eq:expansionformula}
\lambda_1(n) \sim x\sum_{k=0}^{\infty} \frac{C_k(x)}{n^k}\,,
\end{equation}
where $x = j_{0,1}^2$ is the liming eigenvalue (i.e., the lowest eigenvalue of the unit disk, which is given by the square of the first root of the Bessel function~$J_0(x)$). Since regular polygons are approaching the circle in the limit, it is natural that $C_0 = 1$. Interestingly, the remaining coefficients seemed to be expressible in closed-form as polynomials in $x = j_{0,1}^2$ of degree $\lfloor\frac{i-3}{2}\rfloor$ with coefficients expressible in terms of Riemann zeta values
\begin{align*}
C_0(x) &= 1 \,,\\
C_1(x) &= 0 \,,\\
C_2(x) &= 0 \,,\\
C_3(x) &= 4\zeta(3) \,,\\
C_4(x) &= 0 \,,\\
C_5(x) &= \zeta(5)(-2x+12) \,,\\
C_6(x) &= \zeta(3)^2(4x+8) \,,\\
C_7(x) &= \zeta(7)(-\tfrac{1}{2}x^2-12x+36) \,,\\
C_8(x) &= \zeta(3)\zeta(5)(2x^2+8x+48) \,.
\end{align*}
The first coefficients up to third order have been discovered by Grinfeld and Strang~\cite{GRINFELD2012135} in 2012 through deforming the circle to a regular polygon and then applying a technique called \textit{calculus of moving surfaces}. Three years later, Boady in his PhD thesis~\cite{boady} managed to compute two more terms also using methods from the calculus of moving surfaces. The seventh and eighth coefficients were recently found by the second author~\cite{jones2017fundamental} using numerical methods similar to the ones presented in this report.

We managed to obtain eight higher order coefficients in this project
\begin{align*}
C_9(x) &= \zeta(9)(\tfrac{9}{4}x^3+104x^2+438x-1020)
        + \zeta(3)^3(240x+96) \,,\\
C_{10}(x) &= \zeta(7)\zeta(3)(x^3+39x^2-24x+144) 
           + \zeta(5)^2(x^3-6x^2-12x+72) \,,\\
C_{11}(x) &= \zeta(11)(-\tfrac{5}{32}x^4-\tfrac{661}{60}x^3
                       -\tfrac{1623}{20}x^2-176x+372) \\
    &+ \zeta(5)\zeta(3)^2(80x^2+176x+96)
     + \zeta^{\mathrm{sv}}(3,5,3)(\tfrac{1}{5}x^3 + \tfrac{54}{5}x^2)\,,\\
C_{12}(x) &= \zeta(9)\zeta(3)(\tfrac{5}{8}x^4+\tfrac{107}{3}x^3+456x^2
                             -\tfrac{488}{3}x+\tfrac{1360}{3}) \\
          &+ \zeta(7)\zeta(5)(\tfrac{11}{8}x^4+\tfrac{47}{2}x^3
                             -207x^2-216x+432)
           + \zeta(3)^4(-16x^2+\tfrac{272}{3}x+\tfrac{32}{3})\,,\\
C_{13}(x) &= \zeta(13)(-\tfrac{7}{64}x^5-\tfrac{226501}{16800}x^4
                -\tfrac{1283839}{8400}x^3-\tfrac{1447393}{1400}x^2
                -618x + 1260) \\
          &+ \zeta(7)\zeta(3)^2(x^4+34x^3+1236x^2+336x+288) \\
          &+ \zeta(5)^2\zeta(3)(-\tfrac{31}{10}x^4+\tfrac{256}{5}x^3
                                -\tfrac{1128}{5}x^2+336x+288) \\
          &+ \zeta^{\mathrm{sv}}(5,3,5)(-\tfrac{157}{1400}x^4
                                        -\tfrac{549}{350}x^3
                                        -\tfrac{12339}{175}x^2) \\
          &+ \zeta^{\mathrm{sv}}(3,7,3)(-\tfrac{5}{56}x^4
                                        -\tfrac{59}{28}x^3
                                        -\tfrac{747}{14}x^2)\,.
\end{align*}
(The expressions of the coefficients $C_{14}(x)$, $C_{15}(x)$, and $C_{16}(x)$ are given in the Appendix, the numerical expressions as well as the underlying eigenvalue data are provided as an attachment to this paper).
In order to determine these expressions we computed the eigenvalues of 650 $n$-gons to at least 980 digits precision. Through this we also discovered the appearance of single-valued multiple zeta values (MZVs) in the expansion coefficients, starting at eleventh order. In case of $C_{16}(x)$ the basis of single-valued multiple zeta values was found with help of the program \textsc{HyperlogProcedures} developed by Oliver Schnetz. We provide a brief introduction to MZVs and the definitions of the single-valued MZVs that appear in our expressions in the Appendix. We should note that the approach used in this paper does not produce explicitly proven results but rather gives conjectures with very high numerical evidence. We did however manage to prove some of the results for regular polygons with Dirichlet boundary condition in the companion paper \cite{analytics_paper}. Namely, we explicitly derived the results up to (and including) $C_{14}$ and proved that that the expansion coefficients are expressions over the space of multiple zeta values. However, the expressions for $C_{15}$ and $C_{16}$ can only be considered conjectural at this point.

\section{The Method of Particular Solutions}
\setlength{\abovecaptionskip}{-50pt}  
\begin{figure}
\centering
\includegraphics[width=1\textwidth]{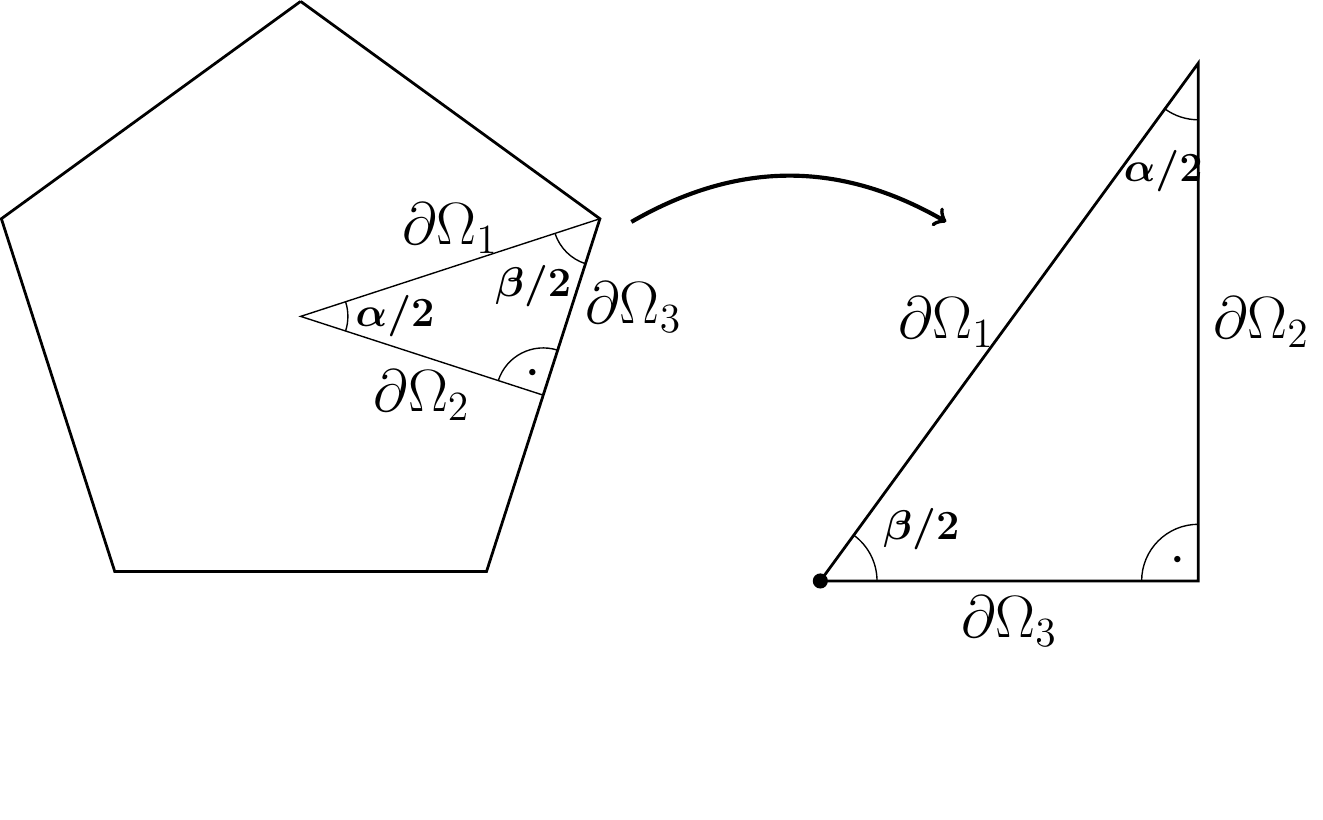}
\caption{Fundamental region of regular polygons for fully symmetric eigenfunctions}
\label{fig:fund_region}
\end{figure}
Because the eigenfunctions considered in this work are dihedrally symmetric, it is sufficient to compute them inside a triangular subdomain (which we refer to as \textit{fundamental domain}) instead of working with the full shape (see Fig.~\ref{fig:fund_region} for an illustration of the fundamental domain for a regular polygon). A numerical method for computing eigenvalues of triangular shapes is given by the method of particular solutions (MPS). The MPS has been introduced by Conway \cite{conway} in 1961 and established by a famous paper of Fox, Henrici and Moler \cite{fhm} in 1966 on computing eigenvalues of L-shape domains.

The main idea of the MPS is to expand $\Psi$ for the spectral parameter $k = \sqrt{\lambda}$ as a Fourier-Bessel series (in polar coordinates)
\begin{equation}
    \Psi(k,r,\theta) = \sum_{\nu=1} c_\nu \psi_\nu(k,r,\theta)\,,
\end{equation}
where
\begin{equation}
\label{eq:ansatz}
\psi_\nu(k,r,\theta) = J_{m_\nu}(k\cdot r)\cdot
\begin{cases}
\sin(m_\nu \theta)\quad &|\quad  \textrm{Odd eigenfunctions}\,,\\
\cos(m_\nu \theta)\quad &|\quad  \textrm{Even eigenfunctions}\,,
\end{cases}
\end{equation}
and $J_{m}$ denotes the Bessel function (see for example \cite{vekua}). The basis functions $\psi_\nu$ have the useful property that they can serve as eigenfunctions along an unbounded wedge by a well-made choice of $m_\nu$. Consider for example a wedge with angle $\alpha$ with the vertex at the origin. Dirichlet boundary conditions impose that the function has to vanish along the boundary. If one chooses $m_\nu = \nu \pi / \alpha$ for the odd eigenfunction basis, one obtains functions that trivially fulfill the boundary condition along both edges of the wedge. This property (combined with the exponential decay of the basis functions) makes the MPS very useful when computing eigenvalues of triangular shapes: the boundary condition on two of the three edges can be trivially fulfilled by the construction of $\psi_\nu$. To obtain the spectral parameter $k$ one truncates the expansion of $\Psi$ to some finite order~$N$ and searches for the parameter $k$ for which the expansion vanishes on a discrete set of points on the remaining edge up to the desired numerical precision (this approach is often referred to as \textit{point-matching}). This results in a linear system of equations
\begin{eqnarray*}
\begin{pmatrix}
\psi_1(k,r_1,\theta_1) & \dots &  \psi_N(k,r_1,\theta_1)\\
\vdots & \ddots & \vdots\\
\psi_1(k,r_N,\theta_N) & \dots & \psi_N(k,r_N,\theta_N)  \\
\end{pmatrix}
\cdot \left( \begin{array}{c}c_1\\ \vdots \\c_N\\\end{array} \right)
= \left( \begin{array}{c}0\\ \vdots \\0\\\end{array} \right)\,,
\end{eqnarray*}
which is square if the number of point-matching points is also chosen to be $N$. The parameter~$k$ now corresponds to an approximation of the true spectral parameter if the coefficients are essentially invariant under the choice of matching points. In view of this, one computes the value of $k$ for which the determinant becomes zero (i.e., when the matrix becomes singular) using root-finding algorithms (we used the secant method) \cite{fhm}. By increasing~$N$, one can obtain better approximations and hence better precision for the eigenvalue.

This approach has been applied by the second author \cite{jones_main,jones2017fundamental} who expanded from the vertex with angle $\beta/2$ (see Fig.~\ref{fig:fund_region}) and applied point-matching to the edge $\partial \Omega_2$  to compute eigenvalues of polygons with 
$n \leq 150$ to 50 digits of precision.
We remark that the MPS has the limitation that the point-matching matrix becomes 
ill-conditioned for large $N$. This has been further analyzed by Betcke and Trefethen \cite{10.2307/20453663, betcke_gsvd} in 2005 who proposed an improved approach that is referred to as \textit{modified MPS}. Their main idea has been to compute an additional matrix that consists of (randomly chosen) points in the interior of the considered region. This matrix is designed to ensure that the eigenfunction is non-zero in the interior. By minimizing the lowest generalized singular value of these two matrices, one can reliably locate the eigenvalues of a given shape because spurious solutions are excluded and $N$ can be chosen as required, without conditioning the matrix too poorly. The modified MPS might be regarded as a successor of the MPS and is superior in most applications. For the computation of eigenvalues of relatively simple shapes at arbitrary precision arithmetic, the original MPS can however still prove itself to be very competitive because the ill-conditioning can be overcome by increasing the working precision and because the generalized singular value decomposition is very computationally expensive compared to the LU-decomposition that is required for the determinant computation. Additionally, the function of the determinant w.r.t. $k$ becomes almost linear close to the eigenvalues which makes the location of the roots very efficient \cite{jones_main}.

We further remark that recent progress has been achieved in making the results of the MPS rigorous (i.e., with certified error bounds) \cite{https://doi.org/10.48550/arxiv.1911.06758,MR4164074}. For our application, certifying the eigenvalues yields no benefit because the derivation of the closed-form expressions using the LLL algorithm is non-rigorous (unless one knows a bound of the height of the coefficients in advance, which we are unaware of). Additionally, specifically certifying the eigenvalues results in a loss of precision.

\section{Domain Decomposition}
Despite the absolute convergence of the MPS (i.e., one always gets better estimates by increasing $N$ as long as the working precision is sufficiently increased) the algorithm which has been applied by the second author in \cite{jones_main,jones2017fundamental} has the disadvantage that the convergence rate (which is the eigenvalue precision per amount of matching points) drastically decreases for polygons with more edges (and hence ``thinner'' fundamental regions). To overcome the decrease of precision for polygons with more edges, we used the technique of domain decomposition. Domain decomposition is by no means new in the context of the MPS. The first application was done by J. Descloux and M. Tolley \cite{Descloux1983AnAA} in 1983 who decomposed domains into several subdomains $\Omega_k$ with eigenfunctions $f_k$ and minimized the quotient of the functionals
\begin{equation}
    \label{eq:functional}
    \sum_{k<l}^{N} \int_{\Gamma_{kl}} \left(|f_k-f_l|^2+|\nabla f_k - \nabla f_l|^2 \right)\,ds\,,
\end{equation}
and
\begin{equation}
    \sum_{k=1}^{N} \int\int_{\Omega_k} |f_k|^2\, dx\, dy\,,
\end{equation}
where $\Gamma_{kl}$ denote straight intersecting boundaries. This method has been further developed by Driscoll \cite{driscoll} in 1997 to compute eigenmodes of the famous isospectral drums of Gordon, Webb and Wolpert \cite{gordon1992hear}. Betcke \cite{10.1093/imanum/drl030} in 2006 formulated the method of Descloux and Tolley as a GSVD problem which does not require the evaluation of the boundary- and domain integrals.

Our approach is based on the idea to treat the intersection of regions as an additional point-matching condition. We chose to decompose the fundamental region of regular polygons into two regions. Then, instead of minimizing Eq.~\eqref{eq:functional}, we explicitly matched the eigenfunctions of both regions and their derivatives on a discrete set of points along the intersecting line. We do not claim this approach to be new (in fact it was already mentioned in Driscoll's paper \cite{driscoll}) however we are unaware of any author proving the success of this method. 

\section{Outline of the Algorithm}
\begin{figure}
\centering
\includegraphics[width=1\textwidth]{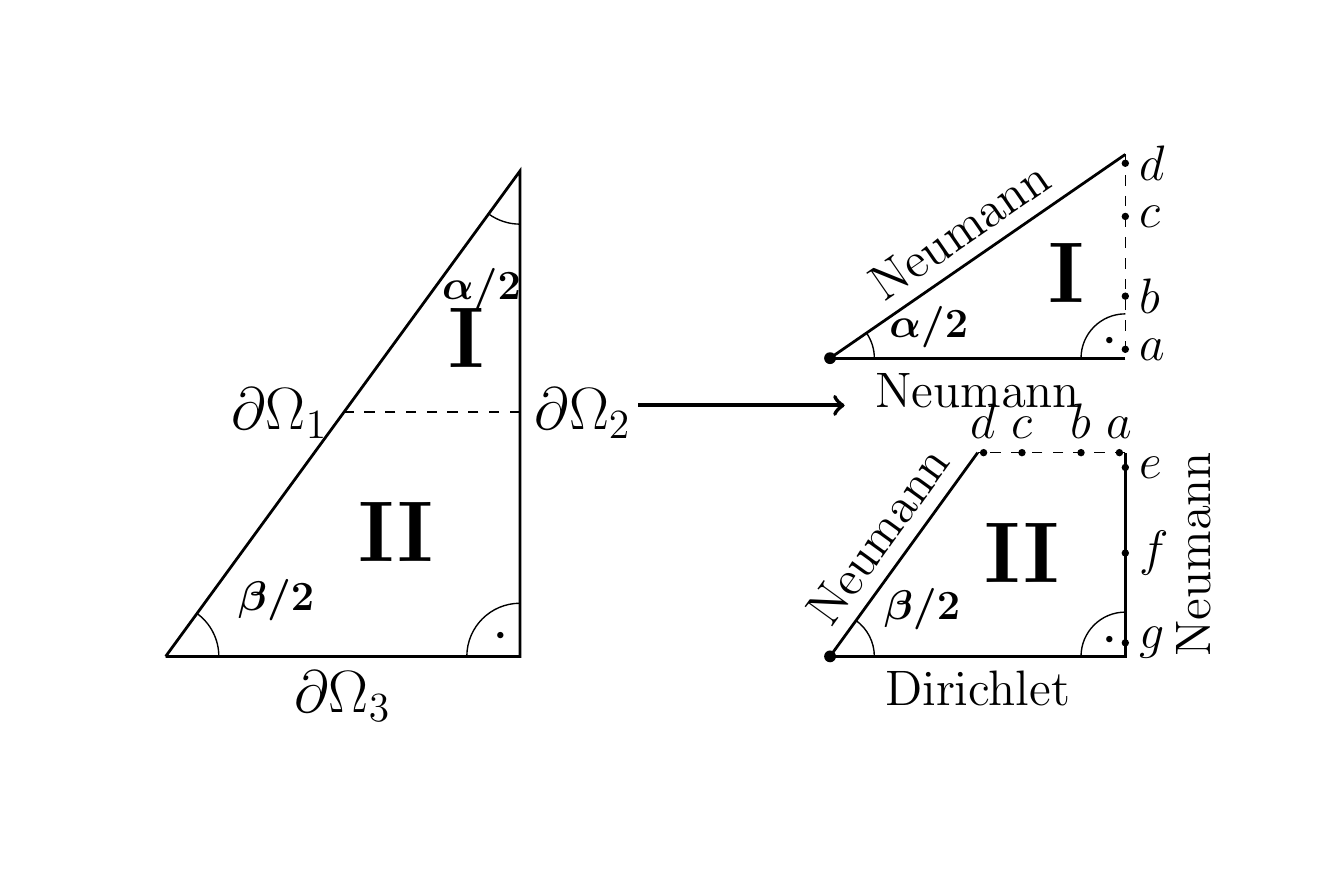}
\caption{Domain decomposition of the fundamental region}
\label{fig:decomposition}
\end{figure}
As shown in Fig.~\ref{fig:decomposition}, the fundamental region can be decomposed into two regions \uproman{1} \& \uproman{2}. For the first region one chooses the eigenfunctions to be
\begin{equation}
    \Psi^{\left[N_{\uproman{1}}\right]}_{\uproman{1}}(k,r,\theta) = \sum\nolimits_{\nu=1}^{N{_{\uproman{1}}}} J_{m_\nu}(k\cdot r)\cdot \cos(m_\nu \theta) \,,
\end{equation}
and expands of the vertex with angle $\alpha/2$. To fulfill the Neumann boundary condition at angle $\theta = 0$ (which is required for fully symmetric eigenmodes), the angular derivative of the wave function (which is a sine with a factor that does not matter here) has to vanish. This is always fulfilled since $\sin(m_\nu\cdot 0) = 0$ for all $m_\nu$. The second Neumann boundary condition is fulfilled if $\sin(m_\nu \frac{\alpha}{2}) = 0$ from which follows
\begin{equation}
    m_\nu = \frac{2\pi}{\alpha}(\nu - 1)\,.
\end{equation}
Similarly for the second domain one expands from the $\beta/2$-vertex (and places it at the origin). We denote this coordinate system with a tilde. The space of eigenfunctions now becomes
\begin{equation}
    \Psi^{\left[N_{\uproman{2}}\right]}_{\uproman{2}}(k,\tilde{r},\tilde{\theta}) = \sum\nolimits_{\mu=1}^{N{_{\uproman{2}}}} J_{m_\mu}(k\cdot \tilde{r})\cdot \sin(m_\mu \tilde{\theta})\,,
\end{equation}
with
\begin{equation}
    m_\mu = \frac{\pi}{\beta}(2\mu - 1)\,.
\end{equation}
The point-matching matrix now has to treat three different boundary conditions:
\begin{enumerate}
    \item The wave functions of both regions have to match along the intersecting line
        \begin{equation*}
        \begin{split}
        \Psi_{\uproman{1}}(a) =&\, \Psi_{\uproman{2}}(a)\,,\\
        \vdots\\
        \Psi_{\uproman{1}}(d) =&\, \Psi_{\uproman{2}}(d)\,.
        \end{split}
        \end{equation*}
    \item The normal derivatives of both wave functions have to match along the intersecting line
        \begin{equation*}
            \begin{split}
            \frac{\partial}{\partial x}\Psi_{\uproman{1}}(a) =& -\frac{\partial}{\partial \tilde{y}}\Psi_{\uproman{2}}(a)\\
            \vdots\\
            \frac{\partial}{\partial x}\Psi_{\uproman{1}}(d) =& -\frac{\partial}{\partial \tilde{y}}\Psi_{\uproman{2}}(d)\,,
            \end{split}
        \end{equation*}
        where the minus sign appears because the expansions have opposite orientation along the intersection line. We computed the derivatives explicitly by differentiating the Bessel and sine functions.
    \item The outer boundary condition of the second region has to be fulfilled
        \begin{equation*}
            \begin{split}
            \frac{\partial}{\partial \tilde{x}}\Psi_{\uproman{2}}(e) =&\, 0\,,\\
            \vdots\\
            \frac{\partial}{\partial \tilde{x}}\Psi_{\uproman{2}}(g) =&\, 0\,.
            \end{split}
        \end{equation*}
\end{enumerate}
All these conditions can then be combined into a point-matching matrix
\begin{equation}
   M = 
\left( \begin{array}{@{}c|c@{}}
   \begin{matrix}
      \Psi_{\uproman{1}}(a)\\
      \vdots\\
     \Psi_{\uproman{1}}(d)\\
   \end{matrix}
      &
     \begin{matrix}
      -\Psi_{\uproman{2}}(a)\\
      \vdots\\
     -\Psi_{\uproman{2}}(d)\\
     \end{matrix}
   \\
   \cmidrule[0.4pt]{1-2}
  \begin{matrix}
      \frac{\partial}{\partial x}\Psi_{\uproman{1}}(a)\\
      \vdots\\
     \frac{\partial}{\partial x}\Psi_{\uproman{1}}(d)\\
   \end{matrix}
      &
     \begin{matrix}
      \frac{\partial}{\partial \tilde{y}}\Psi_{\uproman{2}}(a)\\
      \vdots\\
     \frac{\partial}{\partial \tilde{y}}\Psi_{\uproman{2}}(d)\\
     \end{matrix}
   \\
   \cmidrule[0.4pt]{1-2}
   0 &  \begin{matrix}
      \frac{\partial}{\partial \tilde{x}}\Psi_{\uproman{2}}(e)\\
      \vdots\\
     \frac{\partial}{\partial \tilde{x}}\Psi_{\uproman{2}}(g)\\
     \end{matrix} \\
\end{array} \right)
\end{equation}
We chose the height of the intersecting line to be equal to the length of $\partial \Omega_3$ so that region~\uproman{2} becomes close to a square. We have also experimented with different ways of decomposition (such as choosing the intersection line crosswise through the region, which eliminates the third point-matching condition), but none of them provided better results than the presented approach. In order to overcome the ill-conditioning of the point-matching matrix we performed the computations with about $1.2N$ digits working precision where $N$ is the size of the point-matching matrix (this choice was obtained through trial and error, one could also attempt to choose $N$ based on the asymptotic decay of the series, which does however require a growth condition for $c_\nu$). The points along all boundaries (exterior and interior) were distributed using the Chebyshev distribution. By empirical testing we chose the number of matching points for all three boundary conditions equally to be $N/3$ each while the first eigenfunction is expanded up to $N_{\uproman{1}} = N/4$ and $N_{\uproman{2}} = 3N/4$.

\setlength{\abovecaptionskip}{0pt} 
\begin{figure}
    \centering
    \begin{tikzpicture}
\begin{axis}[
    xlabel={Number of Polygon Edges $n$},
    ylabel={Eigenvalue Precision in Digits},
    xmin=0, xmax=1000,
    ymin=0, ymax=90,
    x tick label style={/pgf/number format/.cd, scaled x ticks = false, set thousands separator={}, fixed},
    xtick={0, 100, 200, 300, 400, 500, 600, 700, 800, 900, 1000},
    ytick={0, 10, 20, 30, 40, 50, 60, 70, 80, 90},
    legend style={at={(1,0.85)},anchor=north east},
    ymajorgrids=false,
    grid style=dashed,
    height=9cm,
	width=12cm,
]

\addplot coordinates {
    (10, 81)(20, 74)(30, 80)(40, 78)(50, 78)(60, 78)(70, 78)(80, 78)(90, 80)(100, 80)(200, 80)(300, 82)(400, 82)(500, 82)(600, 82)(700, 82)(800, 83)(900, 83)(1000, 84)
    };

\addplot coordinates {
    (10, 81)(20, 41)(30, 29)(40, 23)(50, 20)(60, 18)(70, 16)(80, 14)(90, 14)(100, 12)(200, 8)(300, 7)(400, 6)(500, 6)(600, 1)(700, 1)(800, 1)(900, 1)(1000, 1)
    };
    
\legend{With Domain Decomposition,Without Domain Decomposition}
    
\end{axis}
\end{tikzpicture}
    \caption{Comparison of the eigenvalue precision for different polygons with 200 matching points}
   \label{fig:precisions}
\end{figure}
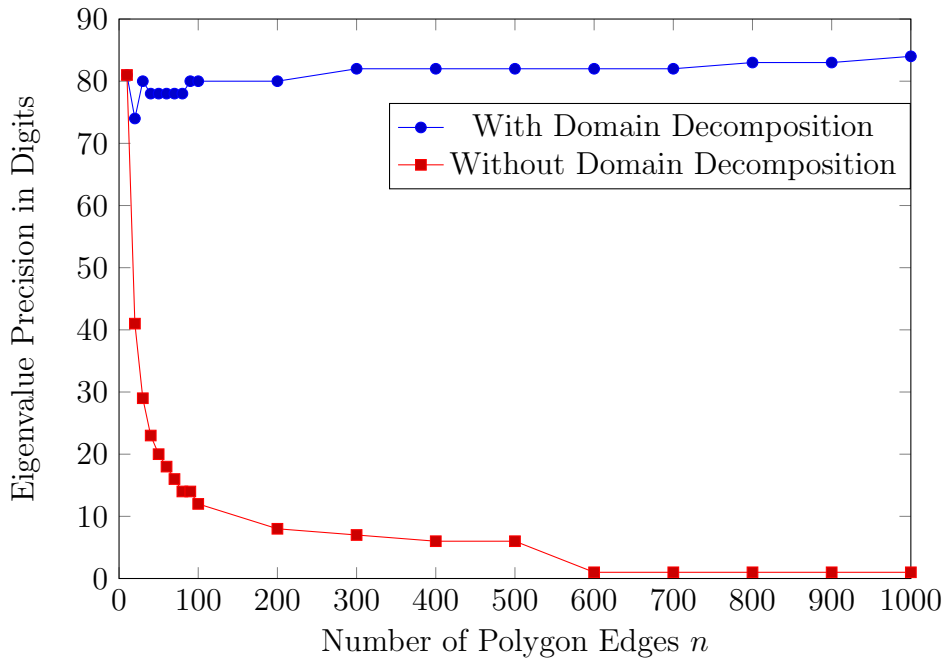
The advantage of this approach becomes clear in Fig.~\ref{fig:precisions}: the algorithm with domain decomposition is almost unaffected by the number of polygon edges $n$ while the original point-matching algorithm without domain decomposition suffers from a heavy decrease of the convergence rate for increasing~$n$. This allowed us to compute eigenvalues for hundreds of regular polygons to high precision.

\section{Details on the numerical implementation}
We implemented the algorithm in the \textsc{Julia} programming language \cite{DBLP:journals/corr/abs-1209-5145}. We made use of \textsc{Arb} \cite{Johansson2017arb} and its \textsc{Julia} interface implemented in \textsc{Nemo} \cite{Fieker_2017}. \textsc{Arb} is a highly optimized C library for arbitrary precision arithmetic. We used \textsc{Arb} to evaluate the occurring special functions. Especially the evaluation of Bessel functions is computationally demanding (in fact, most computational time was spent on computing these functions) which makes an efficient implementation very important (see \cite{arb_bessel} for some benchmarks of \textsc{Arb}'s Bessel implementation). To compute the determinants we used a wrapper to the function \texttt{arb\_mat\_approx\_lu} of \textsc{Arb} which also runs considerably faster than the default Julia determinant algorithm \cite{arb_linear_algebra}. In general, \textsc{Arb} is designed for interval-arithmetic (which means that every float has an error-bound attached so that rounding errors during arithmetic operations are taken into account). However, for our implementation interval-arithmetic is counterproductive because the point-matching matrix is ill-conditioned and empirical tests have revealed that interval-arithmetic might cancel significant digits too pessimistically during the determinant computation (as it ensures that the result is correct up to its displayed precision). It is therefore important to note that \texttt{arb\_mat\_approx\_lu} is one of the few functions of \textsc{Arb} where approximate (non-interval) arithmetic is used. The computing cluster at which the computation were performed consisted of nodes with 2x~\texttt{Intel Xeon E5-2680 v4 @ 2.40GHz} CPUs with 14 (physical) cores each and 128GB of RAM per node. In total, 39 nodes were available. The computations were queued using \textsc{HTCondor} \cite{htcondor} which offers priority-based distribution of hardware resources between multiple users. This means that the hardware was shared with other users which might have increased CPU-times. Each time the determinant of a point-matching matrix has been computed (within the secant-method) the result (and other significant parameters) were stored and the job re-queued. This approach allowed other users to occupy resources as well and the machines to perform required reboots. If however no other user was using the computing cluster, all resources could be used to prevent idle time. The largest computed matrices were of size $3039\times 3039$ with 3647 digits precision per entry. The computation of these matrices (and their determinant) required around 45 hours of wall-clock time running on 9 threads (in our latest version we computed the point-matching matrix in parallel and then computed the determinant using \textsc{Arb}). To store such a matrix and compute its determinant, about 20 GB of RAM was needed. The computations took several months and required 2.7 million thread-hours in total (the amount of physical core-hours without hyper-threading is impossible to determine).

\section{Derivation of the coefficients}
\label{sec:coefficient_derivation}
To compute the coefficients of the $1/n$-expansion, we used Richardson-extrapolation as it is described in \cite{bender}. The computed eigenvalues $\lambda(n)$ for $n$-sided polygons can be expanded as a series
\begin{align*}
&\lambda(n) = C_0 + C_1\cdot n^{-1} + C_2\cdot n^{-2} + \dots + C_N\cdot n^{-N}\\
&\lambda(n+1) = C_0 + C_1\cdot (n+1)^{-1} + C_2\cdot (n+1)^{-2} + \dots + C_N\cdot (n+1)^{-N}\\
&\vdots\\
&\lambda(n+N) = C_0 + C_1\cdot (n+N)^{-1} + C_2\cdot (n+N)^{-2} + \dots + C_N\cdot (n+N)^{-N}
\end{align*}
were $n$ is the lowest computed polygon and N is the total amount of computed polygons. Known coefficients were subtracted to increase the precision. This creates a square system of equations which we solved with LU-decomposition to obtain numerical estimates for the coefficients $C_i$. Using this approach gives higher precision than a standard polynomial fit as it was applied in \cite{jones2017fundamental}. Still, we were unable to extract the coefficients with nearly the same precision as the computed eigenvalues. For example for the seventeenth coefficient (for which we could not find a closed-form solution) we get a precision of approximately 580 digits while the eigenvalues were computed to at least 980 digits. We tried multiple ways to extract the coefficients with higher precision (for example by trying to exploit the fact that the Richardson matrix is a Vandermonde-matrix (see \cite{koev}) or by performing the computations with rational arithmetic to not implement any rounding errors) but none of them yielded in better results.

An empirical observation is that one seems to get the best results by using about $0.65D$ eigenvalues for the Richardson extrapolation, where $D$ denotes the (estimated) precision of the eigenvalues in digits. This is the reason why we have computed 650 polygons to about 1000 digits precision.

After computing the numerical expressions of the coefficients, we used \textsc{Pari} \cite{PARI2} to find closed-form guesses using the LLL-algorithm \cite{LLL}. For regular polygons with Dirichlet boundary condition, we have also made use of the results of \cite{analytics_paper}, which give closed forms for the coefficients of $x^0$ and $x^1$, to cut down on the search space.

\section{Further Examples}
\label{sec:further_examples}

\subsection{Regular Polygons with Neumann Boundary Condition}
The presented algorithm works equally well for regular polygons with Neumann boundary condition (one only has to slightly adapt the function space for the second region). We have computed the lowest Laplacian eigenvalue of 455 polygons with Neumann boundary condition to at least 700 digits precision. This computed data provides strong evidence that there is a $1/n$-expansion 
\begin{equation}
\label{eq:expansionformula_neumann}
\lambda_1(n) \sim j_{1,1}^2\sum_{k=0}^{\infty} \frac{C_k(j_{1,1}^2)}{n^k}
\end{equation}
for the Neumann case with coefficients (here $x = j_{1,1}^2$)
\begin{align*}
C_0(x) &= 1 \,,\\
C_1(x) &= 0 \,,\\
C_2(x) &= 0 \,,\\
C_3(x) &= 0 \,,\\
C_4(x) &= 0 \,,\\
C_5(x) &= -4\zeta(5)x \,,\\
C_6(x) &= 0 \,,\\
C_7(x) &= \zeta(7)(-2x^2-22x) \,,\\
C_8(x) &= 0 \,,\\
C_9(x) &= \zeta(9)(-\tfrac{3}{2}x^3-\tfrac{52}{3}x^2-\tfrac{296}{3}x) \,,\\
C_{10}(x) &= \zeta(5)^2(2x^3+16x^2) \,,\\
C_{11}(x) &= \zeta(11)(-\tfrac{5}{4}x^4-\tfrac{215}{12}x^3-108x^2-418x) \,,\\
C_{12}(x) &= \zeta(7)\zeta(5)(4x^4+46x^3+176x^2) \,,\\
C_{13}(x) &= -\zeta(13)(\tfrac{35}{32}x^5 + \tfrac{4043}{150}x^4 + \tfrac{1699}{12}x^3 + 590x^2 +
        1732x) \\
        &-\tfrac{48}{5}\zeta(5)^2\zeta(3)x^4
         -\tfrac{2}{25}\zeta^{\mathrm{sv}}(5,3,5)x^4 \,,\\
C_{14}(x) &= \zeta(9)\zeta(5)(\tfrac{15}{4}x^5 + \tfrac{826}{9}x^4 + \tfrac{920}{3}x^3 +
        \tfrac{2368}{3}x^2)\\
        &+ \zeta(7)^2(\tfrac{9}{4}x^5 + 18x^4 + \tfrac{385}{2}x^3 + 484x^2) \,,\\
C_{15}(x) &= -\zeta(15)(\tfrac{63}{64}x^6 - \tfrac{543919}{6720}x^5 + \tfrac{16651961}{18900}x^4 + \tfrac{33557}{36}x^3 + \tfrac{15014}{5}x^2 + \tfrac{35494}{5}x) \\
        &-\zeta(5)^3(\tfrac{20}{3}x^5+\tfrac{1288}{15}x^4+64x^3) 
        -(\zeta^{\mathrm{sv}}(5,3,7)+336\zeta(7)\zeta(5)\zeta(3))
        (\tfrac{1}{14}x^5-\tfrac{13}{35}x^4)\,,\\
C_{16}(x) &= \zeta(11)\zeta(5)(\tfrac{7}{2}x^6+\tfrac{16993}{120}x^5+\tfrac{7088}{5}x^4+1714x^3+3344x^2)\\
        &+ \zeta(9)\zeta(7)(\tfrac{9}{2}x^6+\tfrac{6409}{72}x^5+\tfrac{1231}{9}x^4+\tfrac{6836}{3}x^3+\tfrac{13024}{3}x^2)\\
        &+ \zeta(5)^2\zeta(3)^2(12x^5-192x^4)+\zeta^{\mathrm{sv}}(3,5,3)\zeta(5)(\tfrac{6}{5}x^5-\tfrac{96}{5}x^4)\,.
\end{align*}
Our results thus indicate that the appearance of single-valued MZVs is not restricted to the Dirichlet eigenvalues (although the results for Neumann boundary condition have not been formally proven yet).

\subsection{Regular Star Polygons}

\begin{figure}
    \newcounter{subfloat}
    \renewcommand{\thesubfloat}{\alph{subfloat}}
    \newcommand{\image}[2]{%
      \stepcounter{subfloat}%
      \begin{tabular}[t]{@{}c@{}}
      #2 \\
      (\thesubfloat) #1
      \end{tabular}%
    }
    \centering
    \image{\{11,2\}}{%
  \includegraphics[scale=.2]{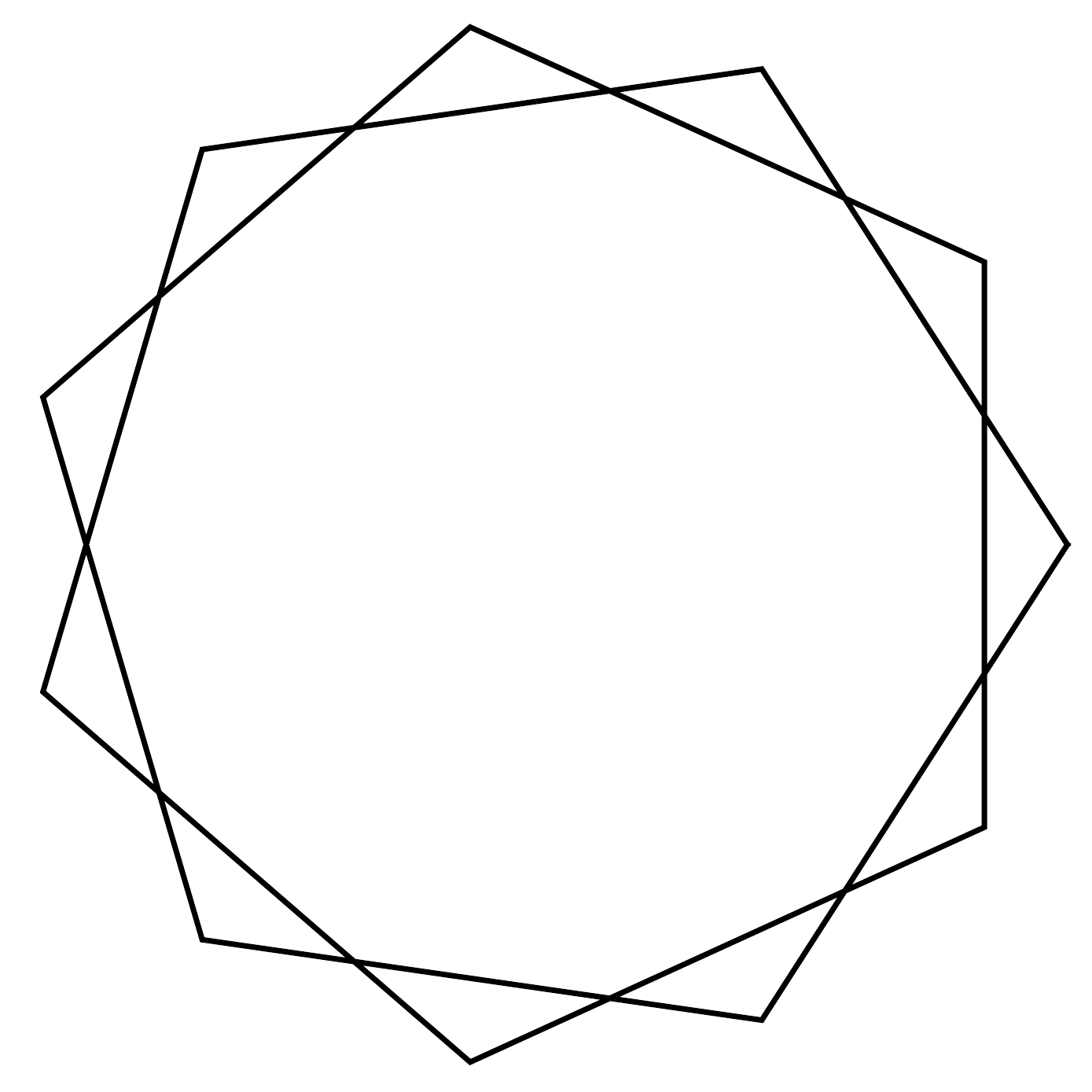}%
    }\quad
    \image{\{11,3\}}{%
      \includegraphics[scale=.2]{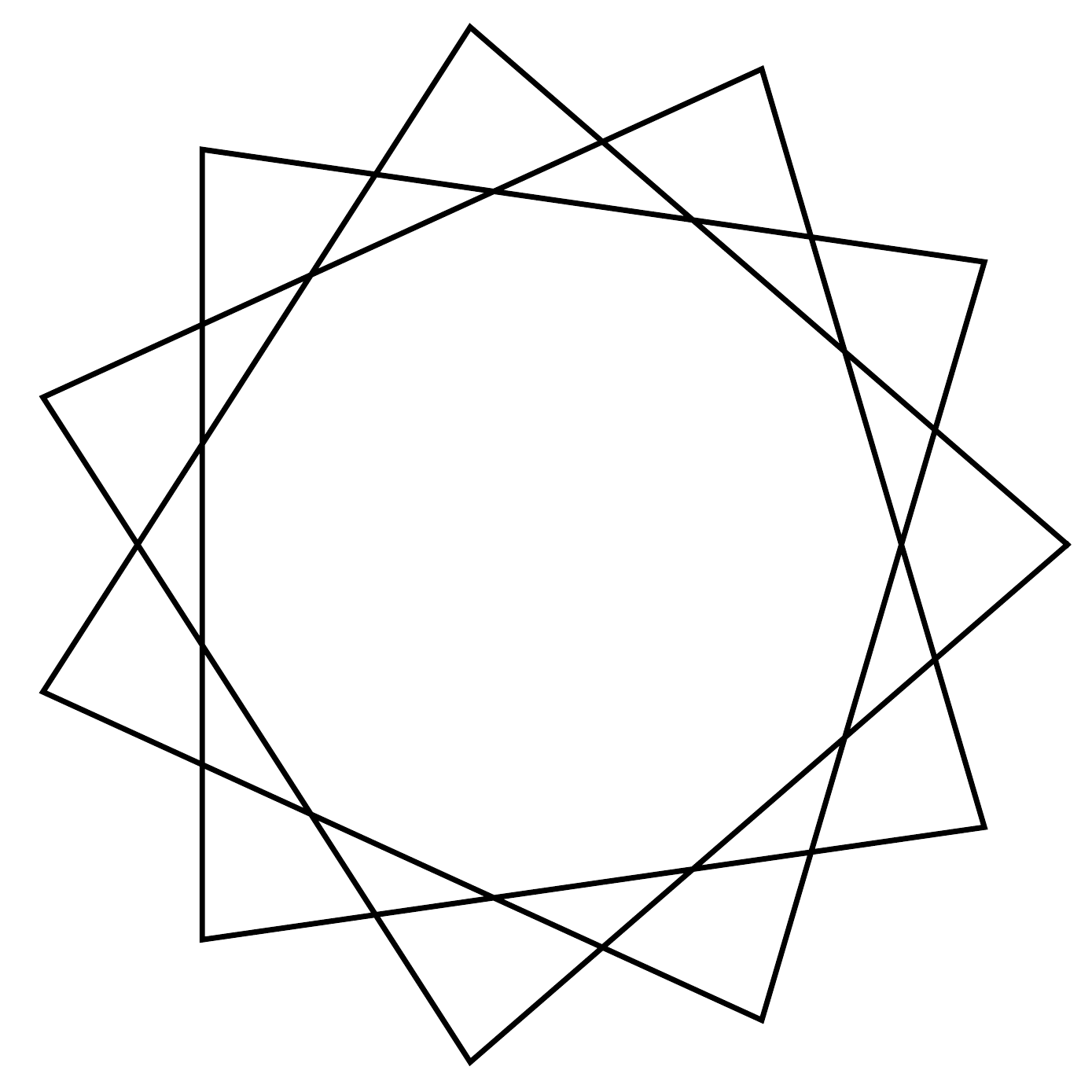}%
    }\quad
    \image{\{11,4\}}{%
      \includegraphics[scale=.2]{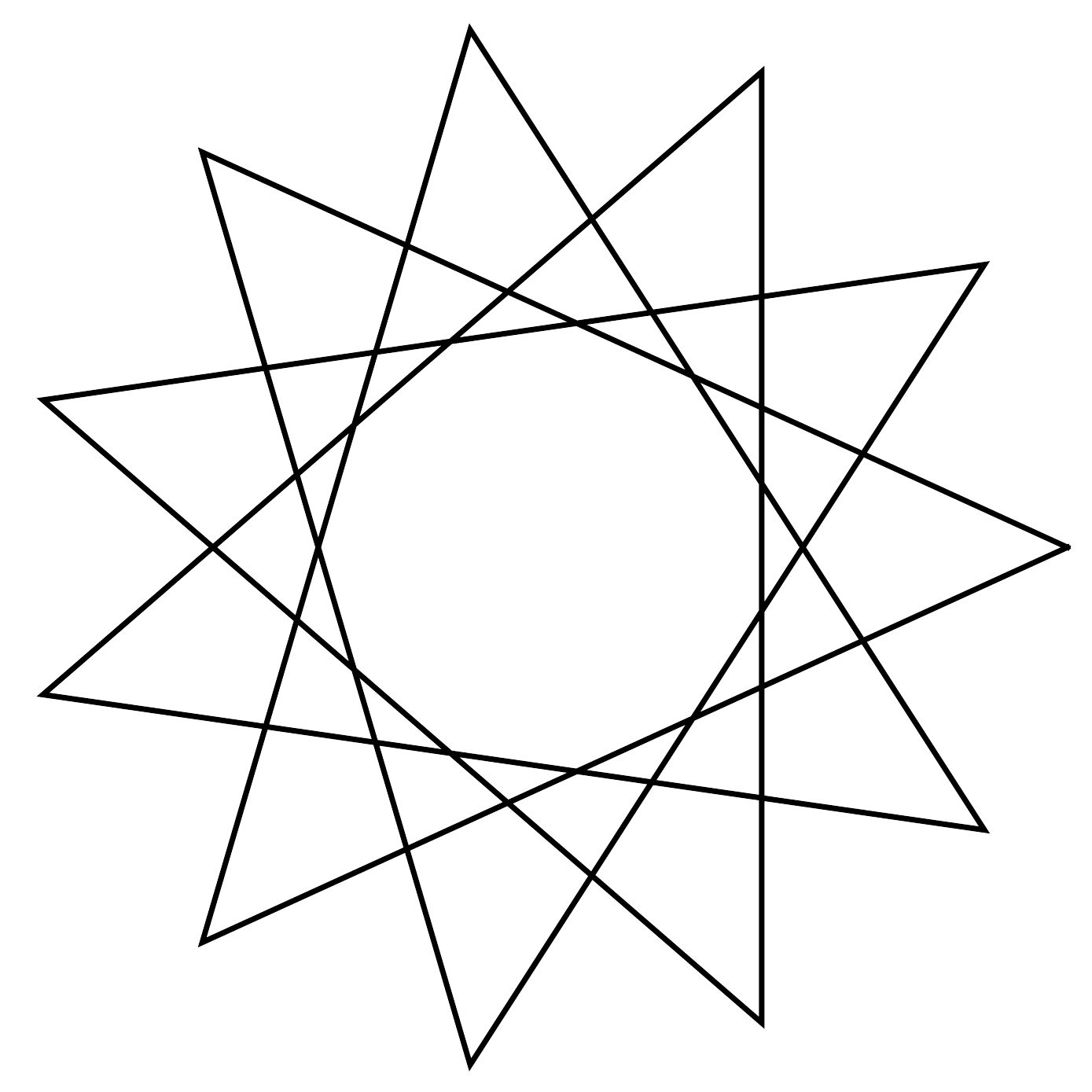}%
    }\quad
    \image{\{11,5\}}{%
      \includegraphics[scale=.2]{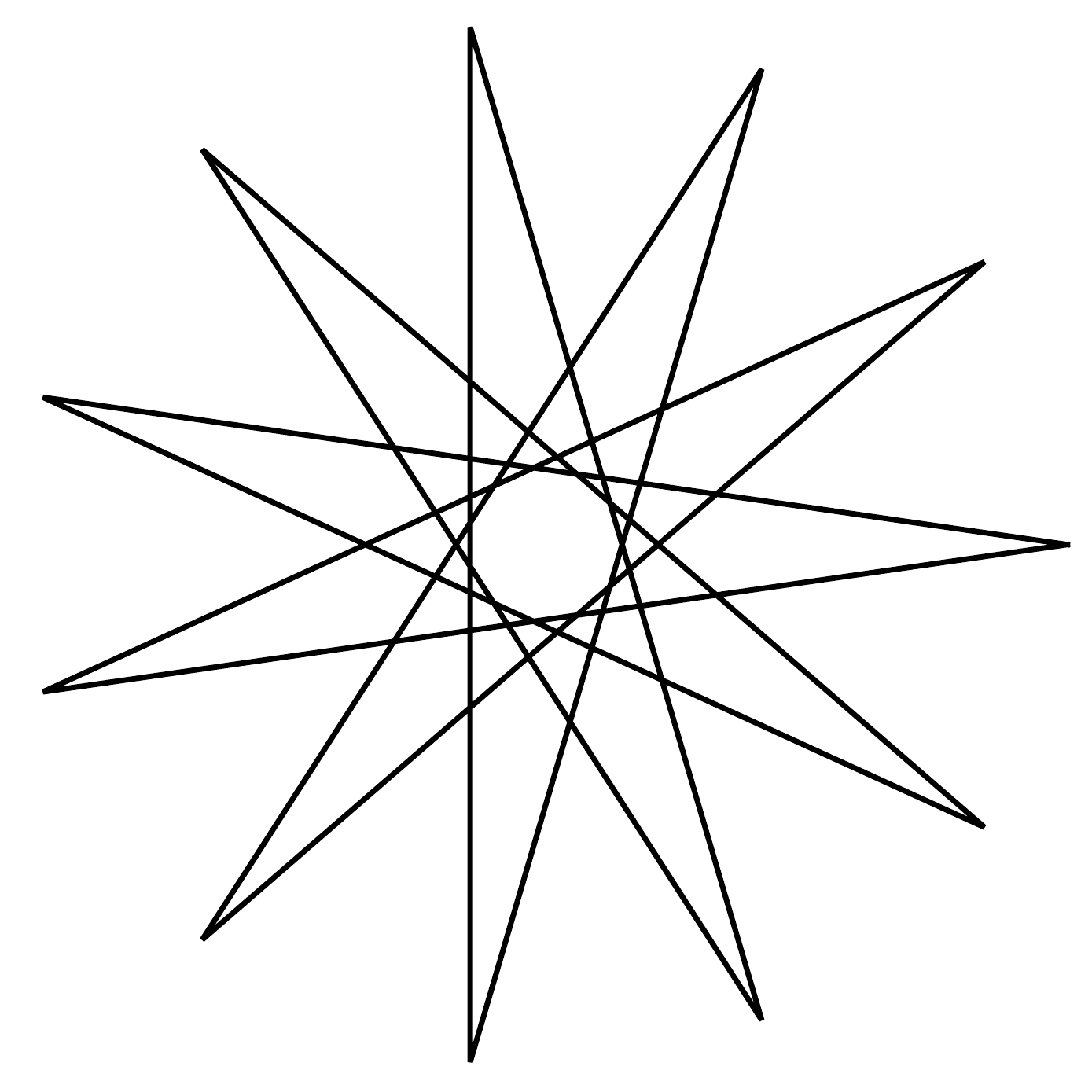}%
    }
    \caption{Regular star polygons with 11 edges and different densities. The figure was created using a \textsc{Mathematica} package by Michael Schreiber \cite{mathematica_star}.}
    \label{fig:star_polygons}
\end{figure}

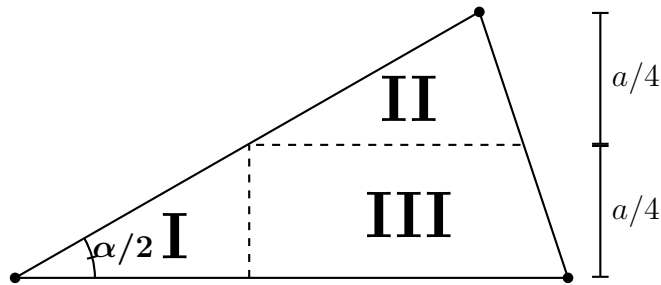
\begin{figure}
	\centering
	\begin{center}
  \begin{tikzpicture}[thick,scale=7]
  \tikzset{shift={(-0.5196426405820805,-0.25190329599000333)}}
  \coordinate (A) at (0,0);
  \coordinate (B) at (1.0392852811641611,0);
  \coordinate (C) at (0.87261861449749445, 0.50380659198000666);

    \draw[arrows=-](A)-- node [above] {} (C);
    \draw[arrows=-](A)-- node [above] {} (B);
    \draw[arrows=-](B)-- node [above] {} (C);
    \draw[arrows=|-|](1.1,0.50380659198000666)-- node [right] { $a/4$} (1.1,0.25190329599000333);
    \draw[arrows=|-|](1.1,0)-- node [right] { $a/4$} (1.1,0.25190329599000333);
    \draw[dashed](0.44,0)-- node [above] {} (0.44,0.25190329599000333);
    \draw[dashed](0.95,0.25190329599000333)-- node [above] {} (0.44,0.25190329599000333);
    \tkzLabelAngle[pos=.21](B,A,C){$\boldsymbol{\alpha/2}$}
    \tkzMarkAngle[size=0.15cm,color=black](B,A,C)
    \fill (A) circle[radius=0.01cm];
    \fill (B) circle[radius=0.01cm];
    \fill (C) circle[radius=0.01cm];
    
    \node[] at (0.3,0.08) {\huge \textbf{I}};
    \node[] at (0.74,0.34) {\huge \textbf{II}};
    \node[] at (0.74,0.12) {\huge \textbf{III}};

  \end{tikzpicture}
\end{center}
	\caption{Domain decomposition of the fundamental region of regular star polygons.}
	\label{fig:decomposition_star}
\end{figure}

In order to study whether the expandability of eigenvalues also applies for other polygonal shapes with dihedral symmetry, we have also considered regular star polygons, which are non-convex regular polygons. For convenient labeling, we make use of the Schl\"afli notation: Any star polygon can be labeled as $\{n,q\}$ where $n$ corresponds to the number of vertices and $q \geq 2$ is referred to as the density (see Fig.~\ref{fig:star_polygons}). The numbers $q$ and $n$ are also required to be relatively prime. The fundamental regions are then given by the triangles with angles $(\tfrac{\pi}{n},\tfrac{\pi(n-2q)}{2n},\tfrac{\pi(n+2q-2)}{2n})$. One example of a regular star polygon is the pentagram which is labeled as $\{5,2\}$. The lowest 79 eigenvalues of this shape have been computed by the second author~\cite{jones_main} to 20 digits of precision. Despite that, we are unaware of any efforts to compute the eigenvalues of these shapes to significantly high precision. The reason for that is probably that they are relatively difficult to compute at high precision since they contain two non-analytic vertices. In fact, we had to compute expansions using three vertices as shown in Fig.~\ref{fig:decomposition_star}. We have computed expansions of the eigenfunctions in each region to roughly equal orders: $N_{\uproman{1}} \approx N_{\uproman{2}} \approx N_{\uproman{3}} \approx N/3$ (we also chose equal numbers of points along each intersection line). The working precision for this triple-expansion approach was chosen to be $0.7N+50$.

We have computed the eigenvalues of star shapes in the case $q=2$ with $12\leq n \leq 263$ to almost 300 digits precision and in the cases $3 \leq q \leq 8$ to about 100 digits of precision. We find that the expansion coefficients of regular star shapes also involve zeta values and can be written as polynomials in $q$. For the first nine coefficients we conjecture that
\begin{align*}
	C_0 &= 1\,,\\
	C_1 &= 0\,,\\
	C_2 &= 0\,,\\
	C_3 &= \zeta(3)(14 q(q-1) + 4)\,,\\
	C_4 &= 0\,,\\
	C_5 &= \zeta(5)(62q(q-1)(q^2-q+1) + 12 - x\cdot\tfrac{31}{4}(q^2 - q + \tfrac{8}{31}))\,,\\
	C_6 &= \zeta(3)^2(2 (7q(q-1) + 2)^2 + x\cdot\tfrac{63}{4}(q^2 - q + \tfrac{16}{63}))\,,\\
	C_7 &= \zeta(7)(254q(q-1)(q^2-q+1)^2+36) - x\cdot A_7q(q-1)(2q-1)^2 \\
        &+ x\cdot\zeta(7)(\tfrac{809}{12}q(q-1)(q^2-q+1)-\tfrac{1571}{16}q(q-1)-12) \\
        &- x^2\cdot\zeta(7)(\tfrac{127}{64}q(q - 1) + \tfrac{1}{2})\,,\\
	C_8 &= \zeta(3)\zeta(5)(28q(q-1)+8)(31q(q-1)(q^2-q+1)+6)\\ 
	    &+ x\cdot \zeta(3)\zeta(5) (110q(q-1)(q^2-q+\tfrac{95}{176})+8)\\
	    &- x\cdot 28\zeta(3)(\zeta(2)\zeta(3)-2\mathrm{Li}_{3,1,1}(-1,1,1))q(q-1)(2q-1)^2\\
	    &+ x^2\cdot \zeta(3)\zeta(5)(\tfrac{255}{32}q(q-1)+2)\,.
\end{align*}
Here the number $A_7$ is defined as
    \[A_7 = \tfrac{56}{9}\mathrm{Li}_{2,2,3}(-1,-1,-1)+8\,\mathrm{Li}_{2,2,3}(-1,1,-1)+\tfrac{7435}{288}\zeta(5,2)+\tfrac{775}{144}\zeta(4,3)\,.\]
Note that these expressions also reproduce the expansion coefficients of regular convex polygons which correspond to $q=1$. What seems particularly interesting is that the coefficients $C_7$ and $C_8$ involve alternating multiple zeta values which did not appear in the case of regular polygons.

\subsection{Smooth Star Shapes}
\begin{figure}
\centering
\begin{center}
  \begin{tikzpicture}[thick,scale=2]
  \coordinate (O) at (0,0);
  \coordinate (A) at (0,3);
  \def\r{1} 
  \def\c{0.0} 

  \draw[black,domain=0:360,samples=500, line width=1] 
       plot ({(1 + 0.3*cos(5*\x))*cos(\x)},{(1 + 0.3*cos(5*\x))*sin(\x)});
  \draw[black,dashed, line width=0.5] (0,0) circle (\r);
  \draw[arrows=|-|](0,0)-- node [above] {$R$} (1,0);
  \draw[arrows=|-|](1,0)-- node [above] {$d$} (1.3,0);

  \end{tikzpicture}
\end{center}
\caption{Smooth Star Shape with $n = 5$ and $d = 0.3$}
\label{fig:curved_star}
\end{figure}
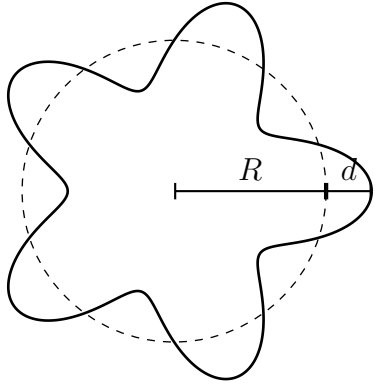

We also investigated the asymptotic expansion of the eigenvalues of shapes with curved vertices which we call smooth star shapes. Smooth star shapes are cycloids with radius 
\begin{equation}
    r(\theta) = R + d \cdot \cos(n\cdot\theta)
\end{equation}
where $n$ corresponds to the number of arcs, $d$ is the height of an arc and $R$ is the radius of the underlying circle (see Fig.~\ref{fig:curved_star}). The eigenvalues of these shapes have previously been investigated by Laura in 1964 \cite{laura} and Wagner in 1971 \cite{wagner} who derived the approximate formula
\begin{equation}
\lambda(n,d) = c^2\cdot \left[ 1 + 0.125~c \cdot\left(J_2(c)/J_1(c) + 2(J_{n-1}(c) - 2J_{n+1}(c))/J_n(c) \right)\cdot d^2 \right]\,,
\end{equation}
with $c = j_{0,1}$. This formula gives an approximation of the eigenvalue at least up to the first digit for $d = 0.3$. Smooth star shapes are relatively easy to construct and have been used to test numerical algorithms for the computation of Laplacian eigenvalues, see for example \cite{guidotti,chen_jiang_chen_yao_2015}. In our case, we investigate the behavior of smooth star shape domains with $d = \left(1/n\right)^m$ where $m\in\mathbb{N}$. This choice of radii ensures that the considered shapes have the eigenvalue of the circle as their limiting eigenvalue. To compute the eigenvalues we simply used an expansion from the interior vertex. This might not be the most efficient approach (since the point-matching matrix becomes extremely ill-conditioned) however it was sufficient to investigate the first terms of the asymptotic expansion. We find that the eigenvalues of these (artificially constructed) shapes can also be written as a $1/n$ series where the expansion coefficients are now rational polynomials of $x = j_{0,1}^2$. For example for $m=2$ one obtains the rational polynomials
\begin{align*}
C_0(x) &= 1 \,,\\
C_1(x) &= 0 \,,\\
C_2(x) &= 0 \,,\\
C_3(x) &= 1 \,,\\
C_4(x) &= 1/2 \,,\\
C_5(x) &= -(2x+1)/4 \,,\\
C_6(x) &= (2x + 3)/4 \,,\\
C_7(x) &= -(12x^2 + 18x - 119)/96\,,\\
C_8(x) &= (48x^2 - 75x - 8)/96 \,,\\
C_9(x) &= -(144x^3 + 3393x^2 - 84x + 167)/2304\,,\\
C_{10}(x) &= (384x^3 + 2499x^2 + 712x + 1203)/768 \,,\\
C_{11}(x) &= -(14400x^4 + 926100x^3 + 2283120x^2 + 525440x - 390979)/368640\,,\\
C_{12}(x) &= (46080x^4 + 933075x^3 + 1041030x^2 - 95940x - 63551)/92160 \,.
\end{align*}
(The remaining coefficients up to $C_{23}$ as well as the coefficients for $m = 3$ and $m=4$ can be found in Appendix.) We omit the case $m=1$ simply for convenience (this case is slightly more difficult to compute to high precision) as we do not expect any substantially different behavior.

\subsection{Hypocycloids}
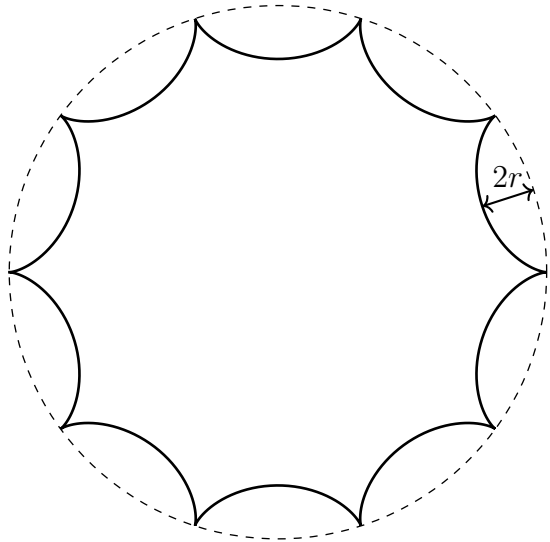
\begin{figure}
	\centering
	\begin{center}
  \begin{tikzpicture}[thick,scale=3]
  \coordinate (O) at (0,0);
  \coordinate (A) at (0,3);
  \def\r{0.1178} 
  \def\n{10} 

  \draw[black,domain=0:360,samples=500, line width=1] 
       plot ({\r*(\n-1)*cos(\x) + \r*cos((\n-1)*\x)},{\r*(\n-1)*sin(\x) - \r*sin((\n-1)*\x)});
  \draw[black,dashed, line width=0.5] (0,0) circle (1.1785);
  \draw[arrows=<->](1.1208309, 0.36418003)-- node [above] {$2r$} (0.8966647, 0.29134402);

  \end{tikzpicture}
\end{center}
	\caption{Hypocycloid with $n=10$}
	\label{fig:hypocycloid}
\end{figure}

\begin{figure}
	\centering
	\begin{center}
  \begin{tikzpicture}[thick,scale=7]
  \coordinate (O) at (0,0);
  \coordinate (A) at (0,0);
  \coordinate (B) at (1.27851, 0);
  \coordinate (C) at (0.8966647, 0.29134402);
  \def\r{0.1178} 
  \def\n{10} 

  \draw[black,domain=0:18,samples=500, line width=1] 
       plot ({\r*(\n-1)*cos(\x) + \r*cos((\n-1)*\x)},{\r*(\n-1)*sin(\x) - \r*sin((\n-1)*\x)});
  \draw[arrows=-](0.8966647, 0.29134402)-- (0, 0);
  \draw[arrows=-](1.17851, 0)-- (0, 0);
  \draw[dashed](0.60531,0)-- node [above] {} (0.60531,0.19667);
  \tkzLabelAngle[pos=.25](B,A,C){$\boldsymbol{\alpha/2}$}
  \tkzMarkAngle[size=0.15cm,color=black](B,A,C)
  \node[] at (0.5,0.08) {\huge \textbf{I}};
  \node[] at (0.8,0.12) {\huge \textbf{II}};
  \fill (0,0) circle[radius=0.01cm];
  \fill (1.17851,0) circle[radius=0.01cm];

  \end{tikzpicture}
\end{center}
	\caption{Decomposition of the fundamental region of a hypocycloid}
	\label{fig:hypocycloid_decomposition}
\end{figure}
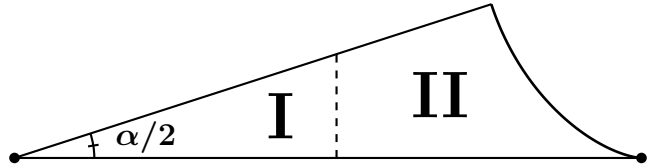

A hypocycloid is the locus of a point on a small circle rotating along the interior of the boundary of a larger circle (see Fig.~\ref{fig:hypocycloid}). In case when the ratio between the radii of the larger and smaller circles is an integer $n$, the resulting hypocycloid has the $n$-sided dihedral symmetry. Interestingly, these shapes are related to group theory: the traces of all matrices in $\mathrm{SU}(n)$ lie inside an $n$-sided hypocycloid as was shown by Kaiser~\cite{kaiser}. Hypocycloids can be parameterized by
\begin{align*}
x(\theta) &= r\cdot(n-1)\cdot \cos(\theta) + r\cdot \cos((n-1)\cdot \theta)\\
y(\theta) &= r\cdot(n-1)\cdot \sin(\theta) - r\cdot \sin((n-1)\cdot \theta)
\end{align*}
where $r = R/n$ is the radius of the smaller (rotating) circle. In the limit as $n\to\infty$ one thus obtains the parametric equations of a circle of radius~$R$. It seems therefore natural to investigate if the eigenvalues of hypocycloids can be expressed as a $1/n$-expansion, similar to what was shown in the previous sections. The vertices of hypocycloids are called \textit{cusps}, i.e., the internal angle at each vertex is zero.
To compute the eigenvalues of a hypocycloid one can compute the expansions from the interior angle and from the cusp and match the two eigenfunctions at the overlap 
(see Fig.~\ref{fig:hypocycloid_decomposition}). Using the Fourier-Bessel expansion ansatz it is not possible to directly fulfill the boundary condition along the curved arcs. The choice of the second eigenfunction (centered at the origin) is therefore the general Fourier-Bessel expansion with integer coefficients $m_\nu$
\begin{equation}
\Psi^{\left[N\right]}_{\uproman{2}}(k,r,\theta) = \sum\nolimits_{m_\nu=1}^{N_{\uproman{2}}{}} J_{m_\nu}(k\cdot r)\cdot \cos(m_\nu \cdot \theta)\,.
\end{equation}
The conditions are now that the two eigenfunctions (including their derivatives) have to match along the intersection line and additionally that the second function has to be zero along the arc. As far as we are aware, so far no other author has applied the method of particular solutions using expansions from a cusp vertex. The MPS (combined with domain decomposition) however works decently and achieves convergence with around 1 digit of eigenvalue precision per 25 additional matching points.

We have computed each hypocycloid eigenvalues from 5 to 100 edges for one week on a single thread, which resulted in around 65 digits precision per eigenvalue. If one fits the eigenvalues on a $1/n$-polynomial-expansion, the expansion coefficients are converging towards constant values which are given by
\begin{align*}
C_0 &= 1\\
C_1 &= 0.372285278452919538950832494889900\\
C_2 &= 1.16087789045807729220949677611\\
C_3 &= 1.8723005009420124517917316
\end{align*}
We have so far been unable to find closed-form expressions for these numbers. Using the symmetry of the expansion formulae for higher fully symmetric eigenvalues (see for example \cite{oikonomou2010casimir}), we found indications that the expansion coefficients $C_i$ are polynomials of degree $\left \lfloor{\frac{|i-1|}{2}}\right \rfloor$ in the limiting eigenvalue $x_m = j_{0,m}^2$. Namely, by computing the first four fully symmetric eigenvalues of several hypocycloids, we found evidence that
\begin{align*}
C_0 &= 1\cdot (x_m)^0\\
C_1 &= 0.372285278452919538950832494889900\cdot(x_m)^0\\
C_2 &= 1.16087789045807729220949677611\cdot(x_m)^0\\
C_3 &= 2.4397482604\cdot(x_m)^0-0.0981202685\cdot(x_m)^1\\
C_4 &= 4.9228013\cdot(x_m)^0-0.088251328\cdot(x_m)^1
\end{align*}
However, it remains unclear whether these numbers can be expressed in closed form.


\section{Conclusion}
We have demonstrated how the method of particular solutions combined with domain decomposition can be efficiently applied to compute Laplacian eigenvalues of several two-dimensional shapes with dihedral symmetry to very high precision. Using this numerical data we have computed two additional eigenvalue expansion coefficients of regular polygons with Dirichlet boundary condition as well as the first 17 coefficients for regular polygons with Neumann boundary condition, and showed that single-valued multiple zeta values appear in the expansions of both of these examples. Moreover, we have also shown that Riemann zeta values and alternating MZVs appear in the eigenvalue expansion coefficients of more general regular star shapes, together with a formula relating the coefficients of star shapes with different density. Additionally, we constructed examples of star shapes with sinusoidal boundary and found that the expansion coefficients are defined over $\mathbb{Q}[x]$ which indicates that the $1/n$ expandability of eigenvalues is not restricted to polygonal shapes only.

\section*{Acknowledgements}
\noindent D.B. would like to thank Oliver Freyermuth for his help in carrying out the computations on a computing cluster, Fredrik Johansson for assistance with Arb, and Plamen Koev for discussions on possible improvements of the approach of Section \ref{sec:coefficient_derivation}. D.R. would like to thank Erik Panzer for help with finding a good basis for single-valued MZVs in weight~16. D.R. would also like to thank Steven Charlton for great help with calculations involving alternating MZVs. D.B. acknowledges financial support from the Bonn-Cologne Graduate School of Physics and Astronomy honors branch.

\newpage
\section{Appendix}

\subsection{Brief Introduction to Multiple Zeta Values}
The Riemann zeta function
\begin{equation}
\zeta(s) = \sum_{n=1}^{\infty} \frac{1}{n^s}
\end{equation}
is known to converge for all $s\in \mathbb{C}$ that satisfy $\textrm{Re}(s)>1$. The Riemann zeta function is arguably the most famous function in all of mathematics due to its tight connection to the distribution of prime numbers. For our treatment we restrict ourselves to the special values of $\zeta(s)$ at natural numbers $s \in \mathbb{N}$, $s>1$, which are sometimes simply called~\textit{zeta values}. For even zeta values, Euler has famously showed that
\begin{equation}
	\zeta(2n) = \frac{(-1)^{n+1}B_{2n}(2\pi)^{2n}}{2(2n)!}\,,
\end{equation}
where $B_n$ are Bernoulli numbers. For odd zeta values much less is known. It has been proven by Ap\'ery in 1979 that $\zeta(3)$ is irrational. His proof has not yet been extended to any other odd zeta values. It is known that at least one of $\zeta(5),\zeta(7),\zeta(9),\zeta(11)$ has to be irrational and that there exists an infinite number of irrational odd zeta values. However, it is still unknown whether all odd zeta values are irrational and whether $\pi, \zeta(3), \zeta(5), \dots, \zeta(2n+1),\dots$ have any algebraic relations among them (altough it is widely believed that there are none). As we have seen in Section~\ref{sec:introduction}, the arguments of products of zeta values that appear in the expansion coefficients of polygon eigenvalues add up to the index of the corresponding asymptotic coefficient. For example, the tenth coefficient for the Dirichlet case contains the products $\zeta(7)\zeta(3)$ and $\zeta(5)^2$, and we have $7+3=2\cdot 5 = 10$. This number is called the weight of a zeta product. Note that the product of two zeta values can be written as
\begin{equation}
\zeta(a)\zeta(b) = \sum_{n=1}^{\infty} \frac{1}{n^a}\cdot \sum_{m=1}^{\infty} \frac{1}{m^b} = \sum_{n,m=1}^{\infty} \frac{1}{n^a m^b}\,.
\end{equation}
We can split the last sum into three terms
\begin{equation}
\sum_{n,m=1}^{\infty} \frac{1}{n^a m^b} = \left(\sum_{0<n<m} + \sum_{0<n=m} + \sum_{0<m<n} \right) \frac{1}{n^a m^b}\,.
\end{equation}
One of the terms is given as a zeta value
\begin{equation}
\sum_{n=m=1}^{\infty}\frac{1}{n^a m^b} = \sum_{n=1}^{\infty}\frac{1}{n^{a+b}} = \zeta(a+b)\,.
\end{equation}
The remaining ``crossing'' terms are the so-called multiple zeta values $\zeta(a,b)$ and $\zeta(b,a)$. More generally, the multiple zeta values (MZVs) are defined as
\begin{equation}
\label{eq:mzv}
\zeta(s_1,\hdots,s_k) = \sum_{0<n_1<\hdots<n_k}\frac{1}{n_1^{s_1}\hdots n_k^{s_k}}\,.
\end{equation}
The decomposition
\begin{equation}
\zeta(a)\zeta(b) = \zeta(a,b) + \zeta(a+b) + \zeta(b,a)
\end{equation}
is sometimes called the Nielsen reflection formula.
Multiple zeta values satisfy a wide variety of interesting relations. For instance, the $\mathbb{Q}$-linear span of MZVs inside $\mathbb{R}$ is closed under products, i.e., forms an algebra. An important subalgebra of MZVs is given by single-valued MZVs that were defined by Brown in~\cite{brown_2014}. These have been found to appear in amplitudes of Feynman diagrams of string theory. The relation between single-valued MZVs and odd zeta values is given by
\begin{equation}
\zeta^\textrm{sv}(2n+1) = 2\zeta(2n+1)\,.
\end{equation}
We will not cover the theory of single-valued MZVs in this brief introduction but instead only give an expression for the single-valued MZVs that appear in our expansion formulas
\begin{align*}
\zeta^\textrm{sv}(3,5,3) &= 2\zeta(3,5,3)-2\zeta(3)\zeta(5,3)-10\zeta(3)^2\zeta(5)\,,\\
\zeta^\textrm{sv}(5,3,5) &= 2\zeta(5,3,5)-22\zeta(5)\zeta(5,3)-120\zeta(5)^2\zeta(3)-10\zeta(5)\zeta(8)\,,\\
\zeta^\textrm{sv}(3,7,3) &= 2\zeta(3,7,3)-2\zeta(3)\zeta(7,3)-28\zeta(3)^2\zeta(7)-24\zeta(5)\zeta(5,3) \\&-144\zeta(5)^2\zeta(3)-12\zeta(5)\zeta(8)\,,\\
\zeta^{\mathrm{sv}}(5,3,7) &= 2\zeta(5, 3, 7) - 28\zeta(3, 5)\zeta(7) - 12\zeta(3, 7)\zeta(5)
-336\zeta(3)\zeta(5)\zeta(7) - 78\zeta(5)^3  - 6\zeta(10)\zeta(5)\\
 &- 28\zeta(8)\zeta(7) - 34\zeta(6)\zeta(9) + 110\zeta(4)\zeta(11) + 1001\zeta(2)\zeta(13) \,,\\
\zeta^{\mathrm{sv}}(3,3,9) &=  2\zeta(3, 3, 9) + 12\zeta(3, 7)\zeta(5) + 30\zeta(3, 5)\zeta(7) + 272\zeta(5)^3 - 27\zeta(3)^2\zeta(9) + 1209\zeta(13)\zeta(2) \\
 & + 252\zeta(11)\zeta(4) + 29\zeta(9)\zeta(6) + 6\zeta(7)\zeta(8) + 318\zeta(3)\zeta(5)\zeta(7) \,,\\
\zeta^{\mathrm{sv}}(1,1,3,4,6) &= 2\zeta(1,1,3,4,6) + 2\zeta(1,1,4,6)\zeta(3) - 12\zeta(3)\zeta(3, 7)\zeta(2) + 4\zeta(3, 3, 7)\zeta(2)\\
&+ 20\zeta(3, 3, 5)\zeta(4) - \tfrac{13133}{56}\zeta(5)^3 -\tfrac{2}{3}\zeta(3)^5 - \tfrac{23}{2}\zeta(3)^2\zeta(5)\zeta(4) + 21\zeta(3, 5)\zeta(5)\zeta(2)\\
&- 16\zeta(3)\zeta(5)^2\zeta(2) - \tfrac{481}{10}\zeta(3, 5)\zeta(7) + \tfrac{349}{45}\zeta(7)\zeta(8)
 - \tfrac{118175}{336}\zeta(5)\zeta(10)
 - \tfrac{2927}{9}\zeta(3)^2\zeta(9) \\
 &- \tfrac{56717}{60}\zeta(13)\zeta(2) + \tfrac{84185}{72}\zeta(9)\zeta(6) - \tfrac{27199}{12}\zeta(11)\zeta(4) - \tfrac{28}{5}\zeta(2)\zeta(3, 5, 5) + \tfrac{58}{9}\zeta(3)\zeta(3, 9) \\&+ \tfrac{3992923}{49752}\zeta(3)\zeta(12)
 + 3\zeta(3)^3\zeta(6) - \tfrac{145}{56}\zeta(5)\zeta(3, 7) - 8\zeta(3)\zeta(3, 5)\zeta(4) - \tfrac{12509}{24}\zeta(3)\zeta(5)\zeta(7) 
 \\&-52\zeta(3)^2\zeta(7)\zeta(2)\,.
\end{align*}
Some of the expressions we give also involve alternating multiple zeta values which are special values of multiple polylogarithms
\begin{equation}
\mathrm{Li}_{m_1,\dots,m_k}(z_1,\dots,z_k) = \sum_{0<n_1<\hdots<n_k} \frac{z_1^{n_1}\dots z_k^{n_k}}{n_1^{m_1}\dots n_k^{m_k}}\,,
\end{equation}
when all $z_i$ are set to $\pm1$.

\subsection*{Remaining Coefficients for Regular Polygons with Dirichlet Boundary Condition}
\begin{align*}
    C_{14}(x) &= \zeta(11)\zeta(3)(1488-664x+\tfrac{32902}{5}x^2+\tfrac{20381}{30}x^3+\tfrac{10169}{240}x^4+\tfrac{7}{16}x^5)\\
&+ \zeta(9)\zeta(5)(1360-\tfrac{2824}{3}x-\tfrac{7306}{3}x^2+\tfrac{1300}{9}x^3+\tfrac{467}{8}x^4+x^5)\\
&+ \zeta(7)^2(648-540x-\tfrac{3627}{2}x^2+\tfrac{483}{4}x^3+\tfrac{175}{16}x^4+\tfrac{9}{16}x^5) \\
&+ \zeta(5)\zeta(3)^3(128+\tfrac{2752}{3}x+\tfrac{3664}{3}x^2+40x^3) \\
&+ \zeta^{\mathrm{sv}}(3,5,3)\zeta(3)x^2\tfrac{1296+28x-x^2}{5}\,,\\
C_{15}(x) &= 
(-\tfrac{96}{11}x^3 - \tfrac{864}{11}x^2)\zeta^{\mathrm{sv}}(1,1,3,4,6)
+ (-\tfrac{7}{72}x^5
- \tfrac{71}{36}x^4
- \tfrac{422}{11}x^3
- \tfrac{6372}{11}x^2)\zeta^{\mathrm{sv}}(3,3,9) \\
&+ (-\tfrac{19}{168}x^5
- \tfrac{2173}{1680}x^4
- \tfrac{13452}{385}x^3
- \tfrac{282042}{385}x^2)\zeta^{\mathrm{sv}}(5,3,7) \\
&+ (-\tfrac{21}{256}x^6
+ \tfrac{2008799}{5040}x^5
+ \tfrac{644330219}{100800}x^4
+ \tfrac{12325006817}{103950}x^3
+ \tfrac{26712690853}{11550}x^2
- \tfrac{10852}{5}x
+ \tfrac{21844}{5}) \zeta(15) \\
&+ (-\tfrac{15}{4}x^5
- \tfrac{2183}{18}x^4
- 4824x^3
- \tfrac{158368}{3}x^2
+ \tfrac{2128}{3}x
+ \tfrac{2720}{3}) \zeta(9)\zeta(3)^2 \\
& + (-\tfrac{133}{12}x^5
+ \tfrac{646}{15}x^4
- \tfrac{274916}{55}x^3
- \tfrac{4524336}{55}x^2
+ 864x + 1728) \zeta(7)\zeta(5)\zeta(3) \\
& + (-3x^5
- \tfrac{511}{10}x^4
- \tfrac{420268}{165}x^3
- \tfrac{2076696}{55}x^2
+ 144x + 288) \zeta(5)^3\\
&+ (-\tfrac{1664}{55}x^3
- \tfrac{62528}{165}x^2
+ \tfrac{3136}{15}x
+ \tfrac{128}{15}) \zeta(3)^5  \,,\\
C_{16}(x) &=(\tfrac{21}{64}\zeta(13)\zeta(3)
+ \tfrac{49}{64}\zeta(11)\zeta(5)
+ \tfrac{29}{32}\zeta(9)\zeta(7)) x^6 \\
&+ (\tfrac{446687}{8400}\zeta(13)\zeta(3)
+ \tfrac{42067}{480}\zeta(11)\zeta(5)
+ \tfrac{15073}{288}\zeta(9)\zeta(7)
- 2\zeta(7)\zeta(3)^3
\\&+ \tfrac{32}{5}\zeta(5)^2\zeta(3)^2
+ (\tfrac{67}{350}\zeta^{\mathrm{sv}}(5,3,5)
- \tfrac{1}{7}\zeta^{\mathrm{sv}}(3,7,3)) \zeta(3)
- \tfrac{1}{20}\zeta^{\mathrm{sv}}(3,5,3)\zeta(5) ) x^5 \\
&+ (
\tfrac{3100857}{2800}\zeta(13)\zeta(3)
+ \tfrac{139433}{80}\zeta(11)\zeta(5)
+ \tfrac{11185}{72}\zeta(9)\zeta(7) 
- 20\zeta(7)\zeta(3)^3
- \tfrac{2584}{5}\zeta(5)^2\zeta(3)^2
\\&- (\tfrac{1139}{350}\zeta^{\mathrm{sv}}(5,3,5)
+ \tfrac{51}{28}\zeta^{\mathrm{sv}}(3,7,3)) \zeta(3)
- \tfrac{47}{10}\zeta^{\mathrm{sv}}(3,5,3)\zeta(5)) x^4 \\
&+ ( 
\tfrac{559169}{60}\zeta(13)\zeta(3)
+ \tfrac{2224}{5}\zeta(11)\zeta(5)
+ \tfrac{3721}{6}\zeta(9)\zeta(7)
+ 988\zeta(7)\zeta(3)^3
- 2808\zeta(5)^2\zeta(3)^2
\\&- (\tfrac{522}{5}\zeta^{\mathrm{sv}}(5,3,5)
- 59\zeta^{\mathrm{sv}}(3,7,3)) \zeta(3)
+ \tfrac{552}{5}\zeta^{\mathrm{sv}}(3,5,3)\zeta(5)) x^3 \\
&+ ( 
\tfrac{16964568}{175}\zeta(3)\zeta(13)
- 26454\zeta(5)\zeta(11)
- 40787\zeta(9)\zeta(7)
+ 18056\zeta(7)\zeta(3)^3
+ \tfrac{150864}{5}\zeta(5)^2\zeta(3)^2
\\&- (\tfrac{214488}{175}\zeta^{\mathrm{sv}}(5,3,5)
- \tfrac{8964}{7}\zeta^{\mathrm{sv}}(3,7,3)) \zeta(3)) x^2 \\
&+ ( 
- 2424\zeta(13)\zeta(3)
- 3480\zeta(11)\zeta(5)
- 4184\zeta(9)\zeta(7)
+ 2240\zeta(7)\zeta(3)^3
+ 3360\zeta(5)^2\zeta(3)^2) x \\
&+(
5040\zeta(13)\zeta(3)
+ 4464\zeta(11)\zeta(5)
+ 4080\zeta(9)\zeta(7)
+ 384\zeta(7)\zeta(3)^3
+ 576\zeta(5)^2\zeta(3)^2)\,.
\end{align*}

\subsection*{Remaining Coefficients for Smooth Star Shapes with $m$ = 2}
\begin{align*}
    C_{13}(x) &= -(1209600x^5 + 169412175x^4 + 1584305000x^3 + 991331400x^2 \\
    &- 67630780x - 22055334)/44236800\,,\\
    C_{14}(x) &= (13271040x^5 + 610174215x^4 + 3119672640x^3 + 1238016540x^2\\
        &+ 7880224x + 62326230)/26542080\,,\\
    C_{15}(x) &= -(457228800x^6 + 118856311875x^5 + 2664517542050x^4 + 8258849383400x^3 \\
        &+ 2125616765440x^2 + 77434447468x - 10887717978)/22295347200\,,\\
    C_{16}(x) &= (33443020800x^6 + 2924868777075x^5 + 38003282808750x^4 + 76583141218500x^3 \\
        &+ 12609060010540x^2 + 18024777715x - 83085435396)/66886041600\,,\\
    C_{17}(x) &= -(1810626048000x^7 + 785618449756875x^6 + 34479052407093000x^5 \\
        &+ 286806810669211000x^4 + 393746334459690800x^3 + 41833381428570655x^2 \\
        &- 382662625041172x - 164371240997550)/112368549888000\,,\\
    C_{18}(x) &= (120394874880000x^7 + 17871670961353125x^6 + 468925800452550000x^5 \\
        &+ 2662356634044423750x^4 + 2565871095884593000x^3 + 179378699791862375x^2 \\
        &- 584375985329753x + 732859733409750)/240789749760000\,,\\
    C_{19}(x) &= -(21184324761600000x^8 + 14232900142095075000x^7 + 1079843748016254200000x^6 \\
        &+ 18704185436989446975000x^5 + 75684321204286594467500x^4 \\
        &+ 52194215245123302758750x^3 + 2424740991307168439892x^2 + 8697289229328086768x \\
        &+ 806072685794305875)/1618107118387200000\,,\\
    C_{20}(x) &= (866843099136000000x^8 + 201696166234942312500x^7 + 9330236455363230750000x^6 \\
        &+ 113623492626641309812500x^5 + 337606293009361365696875x^4 \\
        &+ 168745237201472986724250x^3 + 5201634554633731393011x^2 \\
        &+ 7062388870424537486x - 2948155439034622500)/1733686198272000000\,,\\
    C_{21}(x) &= -(1779483279974400000000x^9 + 1750724028590009405625000x^8 \\
    &+ 210914990771238755890000000x^7 + 6570081043597798567754375000x^6 \\
    &+ 58697712321536371300988890000x^5 + 130785026894321633100311176125x^4 \\
    &+ 47808464146250191120470311872x^3 + 976817476078605785634707984x^2 \\
    &- 709576341339684500163424x - 466588962879561812430000)/163105197533429760000000\,,\\
    C_{22}(x) &= (2140755717626265600000000x^9 + 736244702627870619077343750x^8 \\
        &+ 54826137837085868012812500000x^7 + 1224633170486034514394552343750x^6 \\
        &+ 8269403144636044633016590206250x^5 + 14026490060106825798594234715375x^4 \\
        &+ 3774510130770113767485677126049x^3 + 51231956845133382896804013256x^2 \\
        &- 34284798535733782456620476x \\
        &+ 15568092692374386637462500)/4281511435252531200000000\,,\\
    C_{23}(x) &= -(5869945905985953792000000000x^{10} + 8099489427400822789150312500000x^9 \\
    &+ 1456333964071948513322434125000000x^8 + 74056060297895723139610998046875000x^7 \\
    &+ 1237529454152663240037004207720702500x^6 + 6452892881074851056674758533001245925x^5 \\
    &+ 8424983842296738700856453921666897340x^4 + 1676926979683976178266295659995421360x^3 \\
    &+ 15156930607344875353880119232561056x^2 + 3433373767917417187848282765632x \\
    &+ 1242040349979959360991463800000)/632978650587734212608000000000\,.
\end{align*}

\subsection*{Coefficients for Smooth Star Shapes with $m$ = 3}
\begin{align*}
    C_0(x) &= 1 \,,\\
    C_1(x) &= 0 \,,\\
    C_2(x) &= 0 \,,\\
    C_3(x) &= 0 \,,\\
    C_4(x) &= 0 \,,\\
    C_5(x) &= 1 \,,\\
    C_6(x) &= 1/2 \,,\\
    C_7(x) &= -x/2\,,\\
    C_8(x) &= x/2 \,,\\
    C_9(x) &= -(x^2 + 4x + 2)/8\,,\\
    C_{10}(x) &= (2x^2 + 2x + 3)/4 \,,\\
    C_{11}(x) &= -(x+1)(x^2 + 21x - 18)/16 \,,\\
    C_{12}(x) &= (16x^3 + 104x^2 - 25x + 10)/32 \,,\\
    C_{13}(x) &= -(30x^4 + 1920x^3 + 5547x^2 - 192x - 88)/768\,,\\
    C_{14}(x) &= (384x^4 + 7776x^3 + 11523x^2 - 192x - 304)/768\,,\\
    C_{15}(x) &= -(168x^5 + 23520x^4 + 224715x^3 + 184920x^2 - 224x + 64)/6144\,,\\
    C_{16}(x) &= (3072x^5 + 141312x^4 + 758013x^3 + 367512x^2 + 5696x + 8192)/6144\,.\\
    C_{17}(x) &= (-12096x^6 - 3144960x^5 - 70977933x^4 - 235020744x^3 - 69971328x^2 \\
    &- 927744x + 602368)/589824\,,\\
    C_{18}(x) &= (98304x^6 + 8601600x^5 + 113600385x^4 + 246217560x^3 + 46297920x^2 \\
    &- 15360x + 88832)/196608\,.
\end{align*}

\subsection*{Coefficients for Smooth Star Shapes with $m$ = 4}
\begin{align*}
    C_0(x) &= 1 \,,\\
    C_1(x) &= 0 \,,\\
    C_2(x) &= 0 \,,\\
    C_3(x) &= 0 \,,\\
    C_4(x) &= 0 \,,\\
    C_5(x) &= 0 \,,\\
    C_6(x) &= 0 \,,\\
    C_7(x) &= 1\,,\\
    C_8(x) &= 1/2 \,,\\
    C_9(x) &= -x/2\,,\\
    C_{10}(x) &= x/2 \,,\\
    C_{11}(x) &= -x(x+4)/8 \,,\\
    C_{12}(x) &= x(x+1)/2 \,,\\
    C_{13}(x) &= -(x^3 + 22x^2 + 8x + 4)/16\,,\\
    C_{14}(x) &= (2x^3 + 13x^2 + 2x + 3)/4\,,\\
    C_{15}(x) &= -(5x^4 + 320x^3 + 912x^2 + 24x - 144)/128\,,\\
    C_{16}(x) &= (16x^4 + 324x^3 + 480x^2 - 25x + 10)/32\,,\\
    C_{17}(x) &= -x(7x^4 + 980x^3 + 9360x^2 + 7929x - 64)/256\,,\\
    C_{18}(x) &= x(128x^4 + 5888x^3 + 31584x^2 + 16065x - 64)/256\,,\\
    C_{19}(x) &= (-126x^6 - 32760x^5 - 739344x^4 - 2453067x^3 - 773208x^2 \\
    &+ 1536x + 704)/6144\,,\\
    C_{20}(x) &= (3072x^6 + 268800x^5 + 3550080x^4 + 7730685x^3 + 1545624x^2 \\
    &- 1536x - 2432)/6144\,.
\end{align*}


\bibliographystyle{abbrv}
\bibliography{bib.bib}

\begin{thebibliography}{10}

\bibitem{bender}
C.~M. Bender and S.~A. Orszag.
\newblock Advanced mathematical methods for scientists and engineers -
  asymptotic methods and perturbation theory.
\newblock {\em Springer-Verlag, ISBN 0-387-98931-5}, 1999.

\bibitem{analytics_paper}
D.~Berghaus, B.~Georgiev, H.~Monien, and D.~Radchenko.
\newblock {On Dirichlet eigenvalues of regular polygons}.
\newblock {\em arXiv:2103.01057}, 2021.

\bibitem{10.1093/imanum/drl030}
T.~Betcke.
\newblock {A GSVD formulation of a domain decomposition method forplanar
  eigenvalue problems}.
\newblock {\em IMA Journal of Numerical Analysis}, 27(3):451--478, 2006.

\bibitem{betcke_gsvd}
T.~Betcke.
\newblock The generalized singular value decomposition and the method of
  particular solutions.
\newblock {\em SIAM Journal on Scientific Computing}, 30(3):1278--1295, 2008.

\bibitem{10.2307/20453663}
T.~Betcke and L.~N. Trefethen.
\newblock Reviving the method of particular solutions.
\newblock {\em SIAM Review}, 47(3):469--491, 2005.

\bibitem{DBLP:journals/corr/abs-1209-5145}
J.~Bezanson, S.~Karpinski, V.~B. Shah, and A.~Edelman.
\newblock Julia: {A} fast dynamic language for technical computing.
\newblock {\em CoRR}, abs/1209.5145, 2012.

\bibitem{boady}
M.~Boady.
\newblock Applications of symbolic computation to the calculus of moving
  surfaces.
\newblock {\em Phd thesis, Drexel University}, 2016.

\bibitem{brown_2014}
F.~Brown.
\newblock Single-valued motivic periods and multiple zeta values.
\newblock {\em Forum of Mathematics, Sigma}, 2:e25, 2014.

\bibitem{chen_jiang_chen_yao_2015}
C.~S. Chen, X.~Jiang, W.~Chen, and G.~Yao.
\newblock Fast solution for solving the modified helmholtz equation with the
  method of fundamental solutions.
\newblock {\em Communications in Computational Physics}, 17(3):867–886, 2015.

\bibitem{conway}
H.~D. Conway.
\newblock {The Bending, Buckling, and Flexural Vibration of Simply Supported
  Polygonal Plates by Point-Matching}.
\newblock {\em Journal of Applied Mechanics}, 28(2):288--291, 1961.

\bibitem{MR4164074}
J.~Dahne and B.~Salvy.
\newblock Computation of tight enclosures for {L}aplacian eigenvalues.
\newblock {\em SIAM J. Sci. Comput.}, 42(5):A3210--A3232, 2020.

\bibitem{koev}
J.~Demmel and P.~Koev.
\newblock The accurate and efficient solution of a totally positive generalized
  vandermonde linear system.
\newblock {\em SIAM J. Matrix Anal. Appl.,}, 27(1):142–152, 2005.

\bibitem{Descloux1983AnAA}
J.~Descloux and M.~C. Tolley.
\newblock An accurate algorithm for computing the eigenvalues of a polygonal
  membrane.
\newblock {\em Computer Methods in Applied Mechanics and Engineering},
  39(1):37--53, 1983.

\bibitem{driscoll}
T.~A. Driscoll.
\newblock Eigenmodes of isospectral drums.
\newblock {\em SIAM Rev., 39(1), 1–17}, 1997.

\bibitem{Fieker_2017}
C.~Fieker, W.~Hart, T.~Hofmann, and F.~Johansson.
\newblock Nemo/hecke.
\newblock {\em Proceedings of the 2017 ACM on International Symposium on
  Symbolic and Algebraic Computation - ISSAC ’17}, 2017.

\bibitem{fhm}
L.~Fox, P.~Henrici, and C.~Moler.
\newblock Approximations and bounds for eigenvalues of elliptic operators.
\newblock {\em SIAM J. Numer. Anal.}, 1967.

\bibitem{gordon1992hear}
C.~Gordon, D.~L. Webb, and S.~Wolpert.
\newblock One cannot hear the shape of a drum.
\newblock {\em Bull. Amer. Math. Soc.}, 27:134--138, 1992.

\bibitem{GRINFELD2012135}
P.~Grinfeld and G.~Strang.
\newblock Laplace eigenvalues on regular polygons: A series in 1/n.
\newblock {\em Journal of Mathematical Analysis and Applications}, 385(1):135
  -- 149, 2012.

\bibitem{guidotti}
P.~Guidotti and J.~V. Lambers.
\newblock Eigenvalue characterization and computation for the laplacian on
  general 2-d domains.
\newblock {\em Numerical Functional Analysis and Optimization},
  29(5-6):507--531, 2008.

\bibitem{https://doi.org/10.48550/arxiv.1911.06758}
J.~Gómez-Serrano and G.~Orriols.
\newblock Any three eigenvalues do not determine a triangle, 2019.

\bibitem{Johansson2017arb}
F.~Johansson.
\newblock Arb: efficient arbitrary-precision midpoint-radius interval
  arithmetic.
\newblock {\em IEEE Transactions on Computers}, 66:1281--1292, 2017.

\bibitem{arb_bessel}
F.~Johansson.
\newblock Computing hypergeometric functions rigorously.
\newblock {\em ACM Trans. Math. Softw.}, 45(3), 2019.

\bibitem{arb_linear_algebra}
F.~{Johansson}.
\newblock Faster arbitrary-precision dot product and matrix multiplication.
\newblock In {\em 2019 IEEE 26th Symposium on Computer Arithmetic (ARITH)},
  pages 15--22, 2019.

\bibitem{jones_main}
R.~S. Jones.
\newblock Computing ultra-precise eigenvalues of the laplacian within polygons.
\newblock {\em Adv. Comput. Math.}, 2017.

\bibitem{jones2017fundamental}
R.~S. Jones.
\newblock The fundamental laplacian eigenvalue of the regular polygon with
  dirichlet boundary conditions.
\newblock {\em arXiv:1712.06082}, 2017.

\bibitem{kaiser}
N.~Kaiser.
\newblock Mean eigenvalues for simple, simply connected, compact lie groups.
\newblock {\em Journal of Physics A: Mathematical and General},
  39(49):15287–15298, 2006.

\bibitem{laura}
P.~A. Laura.
\newblock On the determination of the natural frequency of 8 star- shaped
  membrane.
\newblock {\em Journal of the Royal Aeronautical Society}, Pol. 68:274--275,
  1964.

\bibitem{LLL}
A.~Lenstra, H.~Lenstra, and L.~Lovász.
\newblock Factoring polynomials with rational coefficients.
\newblock {\em Math. Ann.}, 1982.

\bibitem{htcondor}
M.~J. Litzkow, M.~Livny, and M.~W. Mutka.
\newblock Condor-a hunter of idle workstations.
\newblock {\em Proceedings of the 8th International Conference of Distributed
  Computing Systems}, pages 104--111, 1988.

\bibitem{Molinari}
L.~Molinari.
\newblock On the ground state of regular polygonal billiards.
\newblock {\em Journal of Physics A: Mathematical and General},
  30(18):6517--6524, 1997.

\bibitem{oikonomou2010casimir}
V.~K. Oikonomou.
\newblock Casimir energy for a regular polygon with dirichlet boundaries.
\newblock {\em arXiv:1012.5376}, 2010.

\bibitem{polya}
G.~Pólya and G.~Szegö.
\newblock Isoperimetric inequalities in mathematical physics.
\newblock {\em Annals of Mathematical Studies}, 1951.

\bibitem{mathematica_star}
M.~Schreiber.
\newblock {Star Polygons}.
\newblock \url{https://demonstrations.wolfram.com/StarPolygons/}, 2011.
\newblock [Online; accessed 23-September-2022].

\bibitem{PARI2}
{The PARI~Group}, Univ. Bordeaux.
\newblock {\em {PARI/GP version \texttt{2.11.2}}}, 2019.
\newblock available from \url{http://pari.math.u-bordeaux.fr/}.

\bibitem{vekua}
I.~N. Vekua.
\newblock {\em New methods for solving elliptic equations}.
\newblock North-Holland Series in Applied Mathematics and Mechanics, Vol. 1.
  North-Holland Publishing Co., Amsterdam; Interscience Publishers John Wiley
  \& Sons, Inc., New York, 1967.

\bibitem{wagner}
H.~{Wagner}.
\newblock {Fundamental Frequency of a Star-Shaped Membrane}.
\newblock {\em Zeitschrift Angewandte Mathematik und Mechanik}, 51(4):325--326,
  1971.

\end{thebibliography}
\end{document}